\documentclass[reqno, 10pt]{amsart}

\usepackage[utf8]{inputenc}
\usepackage[T1]{fontenc}
\usepackage[english]{babel}

\newcommand{\citep}[1]{\cite{#1}}
\newcommand{\citesndintkuss}{\cite[p.~48]{kussthesis}}
\newcommand{\citespinorgritsenko}{\cite[Theorem 4.8]{gritsenkoarithmeticallifting}}
\newcommand{\citespinorneukirch}{\cite[Chapter VII, Theorem 13.2]{neukirchalgzt}}
\newcommand{\citespinorzagier}{\cite[Section 6]{zagierquadraticnumberfields}}

\usepackage[top = 0.8in, bottom = 0.5in, left = 0.8in, right = 0.8in]{geometry}

\usepackage{amssymb, amsfonts}

\usepackage{bbm}

\numberwithin{equation}{section}

\newtheorem{theoremcounter}{theoremcounter}[section]

\newtheorem{corollary}[theoremcounter]{Corollary}
\newtheorem{definition}[theoremcounter]{Definition}

\newtheorem{lemma}[theoremcounter]{Lemma}
\newtheorem{proposition}[theoremcounter]{Proposition}
\newtheorem{theorem}[theoremcounter]{Theorem}

\newcommand{\cal}{\ensuremath{\mathcal}}

\newcommand{\tit}{\itshape}

\newcommand{\thup}{\ensuremath{\textrm{th}}}

\newcommand{\nbd}{\nobreakdash-\hspace{0pt}}

\newcommand{\NN}{\ensuremath{\mathbb{N}}}
\newcommand{\ZZ}{\ensuremath{\mathbb{Z}}}
\newcommand{\QQ}{\ensuremath{\mathbb{Q}}}
\newcommand{\RR}{\ensuremath{\mathbb{R}}}
\newcommand{\CC}{\ensuremath{\mathbb{C}}}

\renewcommand{\Re}{\ensuremath{\mathop{\mathfrak{Re}}}}
\renewcommand{\Im}{\ensuremath{\mathop{\mathfrak{Im}}}}

\newcommand{\isdiv}{\ensuremath{\mid}}
\newcommand{\nisdiv}{\ensuremath{\nmid}}

\newcommand{\GL}[1]{\ensuremath{\mathrm{GL}_{#1}}}
\newcommand{\SL}[1]{\ensuremath{\mathrm{SL}_{#1}}}
\newcommand{\Sp}[1]{\ensuremath{\mathrm{Sp}_{#1}}}

\newcommand{\tr}{\ensuremath{\mathrm{tr}}}

\newcommand{\HS}{\mathbb{H}}

\renewcommand{\frak}{\ensuremath{\mathfrak}}
\renewcommand{\cal}{\ensuremath{\mathcal}}
\newcommand{\bboard}{\ensuremath{\mathbb}}

\newcommand{\cF}{\cal{F}}

\newcommand{\bbZ}{\bboard Z}

\newcommand{\td}{\tilde}

\newcommand{\ov}{\overline}

\begin{document}
\title{The functional equation for the twisted spinor $L$-series of genus $2$}

\author{Aloys Krieg}
\address{Lehrstuhl A f\"ur Mathematik, RWTH Aachen University, Templergraben~55, 52056~Aachen, Germany}
\email{krieg@matha.rwth-aachen.de}

\author{Martin Raum}
\address{Max Planck Institute for Mathematics, Vivatsgasse~7, 53111~Bonn, Germany}
\email{mraum@mpim-bonn.mpg.de}

\maketitle

\begin{abstract}
We prove the functional equation for the twisted spinor $L$-series of a cuspidal, holomorphic Siegel eigenform for the full modular group of genus~$2$.  It follows from a more general functional equation, valid for Rankin convolutions of paramodular cuspforms.  A non-vanishing result for Fourier-Jacabi coefficients of the eigenforms in question is the central pillar of the deduction of the former from the latter functional equation.
\end{abstract}

\section{Introduction}

In this paper, we prove the functional equation of twisted spinor $L$-series attached to cuspidal, holomorphic Siegel eigenforms of genus $2$.

In the theory of automorphic forms, a conjecture, stemming from the Langlands' functoriality conjecture, says that $L$-series attached to automorphic forms for the general symplectic group $\mathrm{GSp}_{n}$ should coincide, up to finitely many local factors, with $L$-series for $\GL{2 n}$.  We say that there is a functorial transfer from $\mathrm{GSp}_{n}$ to $\GL{2n}$.  In a more classical setting, the $L$-series attached to a $\mathrm{GSp}_{n}$ automorphic form is also known as the spinor $L$\nbd series attached to it.  Functorial transfer for holomorphic Siegel modular forms can be proved using converse theorems for $\GL{2 n}$ \citep{CPS96, CPS99}.  They imply that an $L$-series originates in a $\GL{2 n}$ automorphic form, if sufficiently many $\GL{m}$ twists are ``nice'' in a sense made clear in \citep{CPS99}.  Being nice includes the existence of a holomorphic continuation to the whole complex plane and a functional equation satisfied by this analytic continuation.  Recently, it was proved that the spinor $L$\nbd series for degree $2$ Siegel modular forms for the full modular group indeed originate in $\GL{4}$-automorphic $L$-series \citep{PSS11}.  One important ingredient in this proof is the functional equation proved in this paper.

Andrianov initiated the intense study of the spinor $L$-series for genus $2$ Siegel modular forms, using classical methods.  He showed that the twisted and non-twisted spinor $L$-series admit a meromorphic continuation \citep{andrianovquadrforms}.  Under additional assumptions, Andrianov proved that these analytic continuations satisfy the expected functional equation.  More precisely, his proof depends on the non-vanishing of a certain Fourier coefficient of the eigenform in question.  The technique he used was integration against an Eisenstein series coming from a $0$-dimensional cusp.  In \citep{kohnenskoruppacertaindirichlets, kohnenkriegsenguptatwistsofdirichletseries}, Rankin-Selberg integrals for Eisenstein series coming from 1-dimensional cusps were introduced into the study of spinor $L$-series.  The benefit was twofold:  First, Andrianov's technical assumption could be weakened significantly.  Second, the study of twists by Dirichlet characters was much facilitated.  Gritsenko generalized this construction based on paramodular forms \citep{gritsenkoarithmeticallifting}.  The technical assumption that he used is satisfied for all Siegel eigenforms.  Namely, he assumed that either the first or some Fourier-Jacobi coefficient with prime index does not vanish.  The twisted case, though, was still open.  Ku\ss\ \citep{kussthesis} partially closed this gap, but again an assumption on the first Fourier-Jacobi coefficient was needed.  In parallel, efforts were made to generalize the results in \citep{kohnenkriegsenguptatwistsofdirichletseries} to specific types of congruence subgroups~\citep{Ro04}.

We extend both Ku\ss's and Gritsenko's approach to prove the functional equation of all twisted spinor $L$-functions $Z_f^\chi$ attached to cuspidal weight $k$ Siegel eigenforms for the full modular group of genus $2$ (see Definition \ref{def:spinorLdefinition} for a precise definition of $Z_f^\chi$ and its completion $\ZZ_f^\chi$).
\begin{theorem}[Main Theorem]
Let $f \ne 0$ be a cuspidal Siegel eigenform of weight~$k$ for $\Sp{2}(\ZZ)$.  Then the completed spinor $L$-function twisted by a primitive Dirichlet character $\chi$ with conductor $N$ satisfies the functional equation
\begin{gather*}
  \ZZ_f ^{\chi} (s)
= 
  (-1)^k
  \frac{G_{\chi} ^4}{N^2} \,
  \ZZ_f ^{\overline \chi} (2 k - 2 - s)
\text{.}
\end{gather*}
\end{theorem}

The proof proceeds along the following lines: In Section~\ref{sec:paramodulargroup}, we define a normal, non-rational extension of a paramodular congruence subgroup.  Siegel modular forms and their twists are discussed in Section~\ref{sec:siegelmodularforms}.  The extended paramodular group defined in Section~\ref{sec:paramodulargroup} plays an important role in Section~\ref{sec:eisensteinseries}, where we analyze an Eisenstein series $\mathbb{E}^{(p)*}$ attached to it.  The normal extension allows us to deduce a functional equation for $\mathbb{E}^{(p)*}$.  Section~\ref{sec:fstintegralrepresentation} contains the computation of a Rankin-Selberg integral, involving this Eisenstein series and two other Siegel modular forms.  A more complicated variant of $\mathbb{E}^{(p)*}$, that shows up on the right hand side of the Eisenstein series' functional equation, enters in Section~\ref{sec:sndintegralrepresentation}.  This section can be considered as the heart of this work.  A further Rankin-Selberg integral is related explicitly to the Rankin convolution, that already showed up in the preceding section.  The considerations in Section~\ref{sec:eisensteinseries} to \ref{sec:sndintegralrepresentation} are brought together in the final Section~\ref{sec:spinorLseries}.  It also contains an important statement on the non-vanishing of certain Fourier-Jacobi coefficients.

The special case $p = 1$ of the considerations in Section~\ref{sec:siegelmodularforms} was treated in \citep{kohnenkriegsenguptatwistsofdirichletseries}.  The study of the Eisenstein series in Section~\ref{sec:eisensteinseries} was pursued in \citep{kussthesis}~($p = 1$) and \citep{gritsenkoarithmeticallifting}~($N = 1$).  Both integral representations given in Section~\ref{sec:fstintegralrepresentation} and \ref{sec:sndintegralrepresentation} of the Rankin convolution, generalize the work of Ku\ss~\citep{kussthesis}.


\section{A normal extension of the paramodular group}
\label{sec:paramodulargroup}

We write $R^{(m,n)}$ for the module of $m \times n$-matrices with entries in a ring $R$.  The transpose of a matrix $M$ will be denoted by $M^\tr$.  The real symplectic group is
\begin{gather*}
  \Sp{2}(\RR)
:=
  \bigl\{ M \in \RR^{(4,4)} \,:\, M^\tr J M = J \bigr\}
\text{,}
\quad\text{where}\quad
  J
:=
 \left(\begin{array}{cc|cc}
    &   &  1 &   \\
     &  &   & 1 \\
  \hline
   -1  &   &  &   \\
     & -1  &   & 
  \end{array}\right)
\text{.}
\end{gather*}
A typical element of the symplectic group will be written \mbox{$M = \left(\begin{smallmatrix} A & B \\ C & D \end{smallmatrix}\right)$}, where $A, B, C, D \in \RR^{(2, 2)}$ and $A = \left(\begin{smallmatrix}\alpha & a_1 \\ a_2 & \alpha' \end{smallmatrix}\right)$.  We will use an analog notation for the entries of $B$, $C$ and $D$.

By $a \equiv b \; (N)$ we mean that $N \isdiv (a - b)$.  For  $N \in \NN$ and $\kappa \in \NN$ satisfying $\kappa \isdiv N$, the full Siegel modular group \mbox{$\Gamma_2 := \Sp{2}(\ZZ) := \Sp{2}(\RR) \cap \ZZ^{(4,4)}$} has the subgroups
\begin{align*}
   \Gamma_{2,\infty}
&
:= \{M \in \Gamma_2
     \,:\,
     c_2 = d_2 = 0,\, \gamma = 0,\, \delta = 1
   \}
\text{,}
\\
   \Gamma_{2,1}\left(N, N^2/\kappa \right)
&
:= \left\{M \in \Gamma_2
     \,:\,
     c_2 \equiv d_2 \equiv 0 \;(N),\,
     \gamma \equiv 0\;\left(N^2/\kappa \right)
   \right\}
\text{,}\quad\text{and}
\\
   \Gamma_{2,1} ^1 \left(N, N^2/\kappa \right)
&
:= \left\{M \in \Gamma_{2,1}\left(N, N^2/\kappa \right)
     \,:\,
     \delta \equiv \pm 1\; (N)
   \right\}
\text{.}
\end{align*}

For $t \in \NN$, let $\Gamma_2^{(t)}$ be the group generated by $\Gamma_{2,1}(t,t)$ and $M_{1 / t}$, where
\begin{gather}
  M_{\eta}
:= 
  \left(\begin{array}{cc|cc}
   1 &   &   &   \\
     & 1 &   & \eta \\
  \hline
     &   & 1 &   \\
     &   &   & 1
  \end{array}\right)
\in
  \Sp{2}(\RR)
\end{gather}
for any $\eta \in \RR$.  This is the full paramodular group of level $t$.  It contains the subgroup $\Gamma_{2,1} ^{(t)} \big(N, N^2/\kappa\big)$ that is generated by $\Gamma_2(N t, N^2 t / \kappa)$ and $M_{1 / t}$.

Throughout the paper, we will assume that $p \in \NN$ is a prime that does not divide $N$.
For later use, whenever $0 \ne \eta \in \RR$, set
\begin{gather}
  W_{\eta}
:=
  \left(\begin{array}{cc|cc}
      &    & 1 &             \\
      &    &   & 1/\eta \\
  \hline
   -1 &    &   &             \\
      & -\eta &   &
  \end{array}\right)
\text{,}
\qquad
  D_\eta
:=
  \left(\begin{array}{cc|cc}
   1 &   &   &             \\
     & \eta &   &             \\
  \hline
     &   & 1 &             \\
     &   &   & 1/\eta
  \end{array}\right)
\text{,}
\quad\text{and}
\\\nonumber
  P_p
:=
  \left(\begin{array}{cc|cc}
   p \sqrt p           & -p/\sqrt p &  &   \\
   (1 - p) / \sqrt p   & \sqrt p    &  &   \\
   \hline
     &  & \sqrt p   & (p - 1) / \sqrt p \\
     &  & p/\sqrt p & p \sqrt p
  \end{array}\right)
\text{.}
\end{gather}

Recall from \citep{gritsenkoparamodularforms} that $\Gamma_{2} ^{(p) *}$ is a normal extension of the paramodular group $\Gamma_{2} ^{(p)}$, generated by $P_p$.  We define a similar extension of its subgroup $\Gamma_{2,1}^{(p)} \big( N p, N^2 p / \kappa \big)$, given above.  Choose $p^* \in \NN$ such that $p p^* \equiv 1 \; (N)$.  Set
\begin{gather}
  H_p (N)
:=
  \left(
  \begin{matrix}
  U^{- \tr} & \\
   & U
  \end{matrix}
  \right)
  P_p
\text{,}
\quad\text{where}\quad
  U
=
  \left(
  \begin{matrix}
  1 + p^* p & -1 \\
  - p^* p   & 1
  \end{matrix}
  \right)
\text{.}
\end{gather}
Note that $
  \left(
  \begin{smallmatrix}
  U^{- \tr} & \\
   & U
  \end{smallmatrix}
  \right)
\in
  \Gamma_2^{(p)} \big(N, N^2 / \kappa \big)$.  Since $N$ divides the $(4,3)^\thup$ entry of $\sqrt p \, H_p (N)$, we may use $H_p(N)$ to define an extension of that group.  We call the group $\Gamma_{2, 1} ^{(p) *} \big( N p, N^2 p / \kappa \big)$ that is generated by $\Gamma_{2, 1} ^{(p)} \big( N p, N^2 p / \kappa \big)$ and $H_p (N)$ an extended paramodular congruence subgroup.

For convenience, write
\begin{gather}
  i_N
:=
  \big[ \Gamma_{2,1} ^{(p) *} \big( p,p \big) \,:\, \Gamma_{2,1} ^{(p) *} \big( N p, N^2 p \big) \big]
\end{gather}
for the index of $\Gamma_{2,1} ^{(p) *} \big( N p, N^2 p \big)$ in $\Gamma_{2,1} ^{(p) *} \big( p,p \big)$.  We denote the entries of a vector $\lambda$ by $\lambda_i$ for $i \in \NN$.  The last row of any matrix in $\Gamma_{2,1}^{(p)*} \big( N p, N^2 p / \kappa \big)$ has one of the following forms:
\begin{align*}
  \left(
  \begin{matrix}
  N p \lambda_1 & N^2 p / \kappa \lambda_2 & N p \lambda_3 & \lambda_4
  \end{matrix}
  \right)
  \text{,}\quad\text{where }p \nisdiv \lambda_4
\text{;}
\\[6pt]
  \left(
  \begin{matrix}
  N p \lambda_1 & N^2 p / \kappa \lambda_2 & N p \lambda_3 & p \lambda_4
  \end{matrix}
  \right)
  \text{,}\quad\text{where }p \nisdiv \lambda_2
\text{;}
\\[6pt]
  \sqrt p
  \left(
  \begin{matrix}
  N \lambda_1 & N^2 p \kappa \lambda_2 & N \lambda_3 & \lambda_4
  \end{matrix}
  \right)
  \text{,}\quad\text{where }p \nisdiv \lambda_1
\text{;}
\\[6pt]
  \sqrt p
  \left(
  \begin{matrix}
  N p \lambda_1 & N^2 p / \kappa \lambda_2 & N \lambda_3 & \lambda_4
  \end{matrix}
  \right)
  \text{,}\quad\text{where }p \nisdiv \lambda_3
\end{align*}
for some primitive $\lambda \in \ZZ^{4}$ satisfying $\gcd(\lambda_4, N) = 1$.  This provides us with a system of representatives of $\Gamma_{2,\infty} ^{(p)} \backslash \Gamma_{2,1} ^{(p)*} \big( N p,N^2 p / \kappa \big)$.

Given a pair of coprime integers $\lambda \in \ZZ^2$, fix $\left(\begin{smallmatrix}A_\lambda & B_\lambda \\ \lambda_1 & \lambda_2 \end{smallmatrix}\right) \in \SL{2}(\ZZ)$.  Set
\begin{gather}
\label{eq:Mlambdadef}
  M_\lambda
:=
   \left(\begin{array}{cc|cc}
    A_\lambda &   & B_\lambda &   \\
              & 1 &           &   \\
    \hline
    \lambda_1 &   & \lambda_2 &   \\
              &   &           & 1
   \end{array}\right)
\quad\text{and}\quad
  M_{d,\gamma}
:=
  \left(\begin{array}{cc|cc}
   1 & -N d p                      &                      &   \\
     & 1                           &                      &   \\
   \hline
     &                             & q                    &   \\
     & N^2 p \gamma \theta         & N d p                & 1
  \end{array}\right)
\text{.}
\end{gather}
For $\nu, \theta \isdiv N$, a set of representatives of $\Gamma_{2,1} \big( N \nu p, (N \nu)^2 p \big) \backslash \Gamma_{2,1} \big( N p, N^2 p / \theta \big)$ is given by
\begin{gather*}
\Big\{ M_{d,\gamma} M_\lambda \,:\, \gamma \in \ZZ / \theta \nu^2 \ZZ,\, d \isdiv \nu,\, \lambda \in \left\{1,\ldots,\nu/d\right\}^{2} \Big\}
\text{.}
\end{gather*}
It is also a set of representatives for $\Gamma_{2,1}^{(p)} \big( N \nu p, (N \nu)^2 p \big) \backslash \Gamma_{2,1}^{(p)} \big( N p, N^2 p / \theta \big)$.

Any Dirichlet character $\chi$ with period $N$ induces a character on $\Gamma_{2,1} (N, N)$  by means of $\chi(M) = \chi(\delta)$.  Thus it induces a character of  $\Gamma_{2,1} ^{(p)} \big( N p, N^2 p / \kappa \big)$.  Any such character can be extended to $\Gamma_{2,1}^{(p)*} \big( N p, N^2 p / \kappa \big)$ via
\begin{equation*}
  \chi(M)
=
  \begin{cases}
  \chi(\delta)\text{,} & \text{if } M \in \Gamma_{2,1} ^{(p)} \big( N p ,N^2 p/\kappa \big) \text{;} \\
  \sqrt{\overline \chi (p)} \, \chi \big(M H_p (N) ^{-1} \big)\text{,} & \text{otherwise,}
  \end{cases}
\end{equation*}
where $\sqrt{\overline \chi (p)}$ is a fixed root.  If $p \equiv 1 \; (N)$, we use the following notation to distinguish different roots.  The character with $\chi^{+} \big( H_p (N) \big) = 1$ will be denoted by $\chi^{+}$.  Given $k \in \NN$, we write $\chi^{k -}$ for the extension satisfying $\chi^{k -} \big( H_p (N) \big) = (-1)^k$.  In general, we will suppress the superscript $k$, which will be clear from the context.


\section{Twisted Siegel modular forms}
\label{sec:siegelmodularforms}

Denote the Siegel upper half space of genus $2$ by
\begin{gather*}
  \HS_2
:=
  \big\{ Z = X + i Y \in \CC^{(2,2)} \,:\, Z = Z^\tr,\, Y > 0 \big\}
\text{,}
\end{gather*}
where $Y > 0$ means that all eigenvalues of $Y$ are positive.  We write
\begin{gather}
\label{eq:HSelements_notation}
  Z
=
  \left(\begin{matrix}\tau & z \\ z & \tau'\end{matrix}\right)
\text{,}
\qquad
  X
=
  \begin{pmatrix}x & u \\ u & x'\end{pmatrix}
\text{,}
\quad
\text{and}
\quad
  Y
=
  \begin{pmatrix}y & v \\ v & y'\end{pmatrix}
\end{gather}
for a typical element of $\HS_2$ and its real and imaginary part.

The Siegel upper half space is a homogeneous domain for $\Sp{2}(\RR)$ with action
\begin{gather*}
  M \langle Z \rangle
:=
  (A Z + B)\, (C Z + D)^{-1}
\text{.}
\end{gather*}
The slash action of $\Sp{2}(\RR)$ on functions $f : \HS_2 \rightarrow \CC$ is given by
\begin{gather*}
  f \bigl|_k \,M (Z)
:=
  \det(C Z + D)^{-k} f\big( M \langle Z \rangle \big)
\text{.}
\end{gather*}
We say that $f$ is a Siegel modular form of weight $k$ with respect to a group $\Gamma'$ and a character $\chi$ of $\Gamma'$, if $\chi(M)\, f \big|_k \,M = f$ for all $M \in \Gamma'$.  We suppress the subscript $k$, if it is clear by the context.

The vector space of all holomorphic Siegel modular forms of weight $k$ with respect to a character $\chi$ associated to a discrete subgroup $\Gamma'$ of $\Sp{n}(\RR)$ is denoted by $[\Gamma', k, \chi]$.  The space of holomorphic cuspforms is denoted by $[\Gamma', k, \chi]_0$.  For details on holomorphic Siegel modular form, the reader is referred to \citep{freitagsiegelschemodulfunktionen}.

Write $\mathbbm{1}_r$ for the trivial Dirichlet character with period $r \in \NN$.  We will define twists of Siegel modular forms in $\big[\Gamma_{2,1} ^{(p) *} (p, p), k, \mathbbm{1}_1 ^{k-} \big]_0$.

Let $\chi$ be a character of period $N$, and assume that $p \equiv 1\; (N)$.  Write 
\begin{gather*}
f(Z) = \sum_{m = 1} ^\infty f_m (\tau, z)\, e^{2 \pi i\, m \tau'}.
\end{gather*}
for the Fourier-Jacobi series of $f \in [\Gamma_2 ^{(p) *}, k, \mathbbm{1}_1 ^{k-}]_0$.  The twist $f_\chi$ of $f$ by $\chi$ is
\begin{gather}
\label{siegelmodularformtwistdef}
  f_\chi(Z)
:=
  \sum_{m = 1}^\infty \chi(m) f_{m}(\tau, z)\, e^{2 \pi i\, m \tau'}
\text{.}
\end{gather}

To simplify the notation write $e(x) := e^{2 \pi i x}$.  A subscript $\mu (N)$ of a sum denotes summation over a set of representatives modulo $N$.  A direct calculation, analogous to the considerations in \citep{kohnenkriegsenguptatwistsofdirichletseries}, yields
\begin{gather}
\label{siegelmodularformtwistproperties}
  f_\chi
=
  \frac{1}{N}
  \sum_{\nu,\, \mu (N)}
  \chi(\nu) e\big( -\nu \mu/N \big)
  f \big| M_{\mu/N}
\text{.}
\end{gather}
From this equality, we conclude that $f_\chi$ is a holomorphic Siegel cuspform in $\big[\Gamma_{2,1} ^{(p) *} \big( N p, N^2 p\big),\,k,\,\big(\chi^2\big)^{k-}\big]_0$.  Indeed, notice that
\begin{gather*}
  H_p(N)^{-1} M_{\mu/N} H_p (N) M_{-\mu (p - p^* y)/N}
\in
  \Gamma_{2,1} ^{(p)} \big( N p, N^2 p \big)
\text{.}
\end{gather*}
From this, we obtain
\begin{align*}
  f_\chi \big | H_p (N)
&
=
  \frac{1}{N}
  \sum_{\nu,\, \mu (N)}
  \chi(\nu) \, e\big( -\nu \mu/N \big) \,
  f \big| M_{\mu/N} H_p (N)
\\[6pt] &
=
  (-1)^k
  \frac{1}{N}\!
  \sum_{\nu,\, \mu (N)}\!\!
  \chi(\nu) \, e\big( -\nu \mu/N \big) \cdot
\\ &
\hspace{6em}
  \cdot f \big| H_p(N)^{-1} \, M_{\mu/N} H_p (N) \,
          M_{-\mu (p - p^* y)/N} \, M_{\mu (p - p^* y)/N}
\\[6pt] &
=
  (-1)^k f_\chi
\text{.}
\end{align*}

If $\chi$ is primitive, a refined representation for $f_\chi$ is given by
\begin{gather}
\label{siegelmodularformtwistwithprimitivechi}
  f_\chi
= \frac{1}{G_{\overline \chi}}
  \sum_{\mu (N)}
  \overline \chi (\mu) \,
  f \big| M_{\mu/N}
\text{,}
\end{gather}
where $G_{\ov \chi} := \sum_{\nu (N)} {\ov \chi}(\nu) \, e\big( \nu / N \big)$ denotes the Gau\ss\ sum associated to $\overline \chi$.

Given $\nu \in \NN$,
\begin{gather}
\nonumber
  f_\chi \big| W_{N \nu p }
\in
  \big[ \Gamma_{2,1}^{(p) *} \big( N \nu p, (N \nu)^2 p \big),\, k,\, 0,\, \big( \overline\chi ^2 \big)^{k-} \big]_0
\quad\text{and}
\\
\label{siegelmodularformtwistWNtaut}
  f_\chi \big| W_{N p}
= 
  \big( G_\chi^2 /N \big) f_{\overline \chi}
\end{gather}
is an immediate consequence of \eqref{siegelmodularformtwistwithprimitivechi}.

In Section \ref{sec:eisensteinseries}, we will need the next vanishing result.
\begin{lemma}
\label{la:vanishingundertranslations}
If $q \isdiv N$ is a prime, then
\begin{gather}
  \sum_{\nu(q)}
  f_\chi \left| M_{(N p)^2 \nu/q} ^\tr \right.
= 0
\text{.}
\end{gather}
\end{lemma}
\begin{proof}
The sum \mbox{$\sum_{\nu (q)} \chi(\mu + N \nu/q)$} vanishes for any $\mu \in \ZZ$.  The equality $W_{Np} \, M^\tr_{(N p)^2 \nu / p} \, W^{-1}_{N p} = M_{-\nu / q}$ is verified easily.  Combining these facts with \eqref{siegelmodularformtwistwithprimitivechi} and \eqref{siegelmodularformtwistWNtaut}, we complete the proof.
\end{proof}


\section{Paramodular Eisenstein series}
\label{sec:eisensteinseries}

In this section we will define an Eisenstein series for the extended paramodular group $\Gamma_{2,1}^{(p)*}(N p, N^2 p / \kappa)$.  It admits a meromorphic continuation (see Lemma \ref{la:eisensteinreformulations}), and satisfies the functional equation given in Corollary \ref{compeisensteinfcteq}.

We will write $Z_1$ for the upper left entry of any $2 \times 2$ matrix $Z$.  Recall that we assume that $\kappa \isdiv N$.  Given $Z \in \HS_2$ and $s \in \CC$ satisfying $\Re(s) > 2$, define the Klingen-Eisenstein series for $\Gamma_{2,1}^{(p) *}\left(N p,  N^2 p / \kappa\right)$:
\begin{gather}
\label{eisensteinseriesdef}
   E^{(p) *} _{N p, N^2 p / \kappa, \chi} (Z, s) 
:= \sum_{M : \Gamma_{2,1} ^{(p)} \backslash \Gamma_{2,1} ^{(p) *} \left(N p, N^2 p / \kappa\right)}
    \hspace{-1em}
    \chi^+ (M) \left(\frac{\det \big(\Im (M\langle Z \rangle) \big)}
                       {\Im (M\langle Z \rangle)_1} \right)^s
\text{.}
\end{gather}
To define a completion of this Eisenstein series, write $L(s, \chi) := \sum_{n = 1}^\infty \chi(n)\, n^{-s}$ for the Dirichlet $L$-series attached to $\chi$.
\begin{gather}
  \mathbb{E}^{(p) *} _{N p, N^2 p / \kappa, \chi} (Z, s) 
:=
  \big( N^2 / \pi \big)^s \, p^{3s/2} (1 + p^{-s}) \, \Gamma(s) L(2s, \chi) \,
  E^{(p) *} _{N p, N^2 p / \kappa , \chi} (Z, s)
\text{.}
\end{gather}
Our proof of its functional equation will be based on the Epstein $\zeta$-function.  Given two matrices or a matrix and a vector, $M$ and $v$, of compatible sizes, define $M[v] := v^\tr M v$.  Using the abbreviation $I_2 := \left(\begin{smallmatrix} 1 & 0 \\ 0 & 1 \end{smallmatrix} \right)$, for \mbox{$Z = X + i Y \in \HS_2$}, set
\begin{align}
\label{eq:PZdefinition}
   P_Z
&
:= \left(\begin{array}{cc}
    Y & 0 \\
    0 & Y^{-1}
   \end{array}\right)
   \left[\left(\begin{array}{cc}
    I_2 & 0 \\
    X & I_2
   \end{array}\right)\right]\text{.}
\end{align}
We will use the relation $P_Z[M^\tr] = P_{M \langle Z \rangle}$, which holds for all $M \in \Sp{2}(\RR)$, at several occasions.

The generalized Epstein $\zeta$-function and its completion associated to a positive definite matrix \mbox{$P \in \RR^{(4, 4)}$} and \mbox{$u, v \in \RR^{4}$} are defined as
\begin{align*}
  \zeta(s, u, v, P)
&:=
  \sum_{\substack{\lambda \in \ZZ^{2n}\\[1pt]
         \lambda + v \ne 0}}
  e\big( u^\tr \lambda \big) \, P[\lambda + v]^{-s}
\qquad\text{and}
\\
  \zeta^* (s, u, v, P)
&:=
  \pi^{-s} \Gamma(s)\, \zeta(s, u, v, P)
\text{.}
\end{align*}
By the investigations in \citep{terrasharmonicanalysis}, the above series representation converges locally absolutely for $\Re(s) > 2$.

\begin{lemma}
\label{la:eisensteinreformulations}
The completed Eisenstein series has a meromorphic continuation in~$s$. Whenever both sides are defined, the following equalities hold:
\begin{align}
\label{eq:firsteisensteinequality}
&
   \big(\pi / N^2 \big)^s \, p^{-3s/2} \, \Gamma(s)^{-1} \, \mathbb{E}^{(p) *} _{N p, N^2 p / \kappa, \chi} (Z, s)
\\[6pt]\nonumber
&\qquad=
  \sum_{\lambda \in \ZZ^4 \setminus \{0\}}
  \hspace{-0.5em}
  \chi(\lambda_4)
  \Bigl(
  p^{-s}
  P_Z \Bigl[\left(\begin{matrix}
       N \lambda_1 \\ N^2 p / \kappa \lambda_2 \\ N \lambda_3 \\ \lambda_4
       \end{matrix}\right)\Bigr]^{-s}
  +
  P_Z \Bigl[\left(\begin{matrix}
       N p \lambda_1 \\ N^2 p / \kappa \lambda_2 \\ N p \lambda_3 \\ \lambda_4
       \end{matrix}\right)\Bigr]^{-s}
  \Bigr)
\end{align}
and
\begin{align}
\nonumber
& \quad
  N^{2 s} \, \mathbb{E}^{(p) *} _{N p, N^2 p / \kappa, \chi} (Z, s)
\\[6pt]
\label{eq:secondeisensteinequality}
&=
  \sum_{\substack{
        \alpha (N),\, \beta (\kappa) \\[1pt]
        \gamma (N),\, \delta (N^2)
       }}
  \hspace{-1em}
  \chi(\delta)
  \Bigl(
  p^{-s/2}
  \sum_{g (p)}
  \zeta^*
        \bigl(
        s,\, 0,\,
        (\frac{\alpha}{N},\, \frac{\beta}{\kappa},\,
         \frac{\gamma}{N},\, \frac{p \delta + N^2 g}{N^2 p})^\tr,\,
        P_Z
       \bigr)
\\[4pt]\nonumber
& \hspace{5em}
  +
  p^{(s/2) - 1}
  \sum_{h (p)}
  \zeta^*
        \bigl(
        s,\, (0,\,\frac{\kappa h}{p},\, 0,\, 0)^\tr,\,
        (\frac{\alpha}{N},\, \frac{p \beta}{\kappa},\,
              \frac{\gamma}{N},\, \frac{\delta}{N^2})^\tr,\,
        P_Z
       \bigr)
  \Bigr)
\\[6pt]
\label{eq:thirdeisensteinequality}
&=
  \sum_{\substack{
        \alpha (N),\, \beta (\kappa) \\[1pt]
        \gamma (N),\, \delta (N^2)
       }}
  \hspace{-1em}
  \chi(\delta)
  \Bigl(
  p^{-s/2}
  \sum_{g (p)}
  \zeta^*
        \bigl(
        s,\, 0,\,
        (\frac{\alpha}{N},\, \frac{\beta}{\kappa},\,
         \frac{\gamma}{N},\, \frac{p \delta + N^2 g}{N^2 p})^\tr,\,
        P_Z
       \bigr)
\\[4pt]\nonumber
& \hspace{3em}
  + 
  p^{- 3s/2}
  \hspace{-0.7em}
  \sum_{\substack{
        h_1,\,h_2,\,h_3 (p)
       }} \hspace*{-1ex}
  \zeta^*
        \bigl(
        s,\, 0,\,
        (\frac{p \alpha + N h_1}{N p},\, \frac{\beta}{\kappa},\,
         \frac{p \gamma + N h_2}{N p},\, \frac{p \delta + N^2 h_3}{N^2 p})^\tr,\,
        P_Z
       \bigr)
 \Bigr)
\text{.}
\end{align}
If $\chi$ is not the trivial character, the meromorphic continuation of $\mathbb{E}^{(p) *} _{N p, N^2 p / \kappa, \chi} (Z, s)$ is holomorphic in $s$.  Otherwise, it is holomorphic except for a simple pole at $s = 2$ with residuum $2 \kappa\, \varphi(N)/N$.
\end{lemma}
\begin{proof}
For the time being, suppose that the series occurring in the statement are absolutely convergent.  Then
\begin{align*}
& \quad
  (1 + p^{-s}) \, L(2s, \chi) \,
  E^{(p) *} _{N p, N^2 p / \kappa, \chi} (Z, s)
\\[6pt]
& =
  p^{-s}
  \hspace{-0.5em}
  \sum_{\lambda \in \ZZ^4 \setminus \{0\}}
  \hspace{-0.5em}
  \chi(\lambda_4)
  P_Z \Bigl[\left(\begin{matrix}
       N \lambda_1 \\ N^2 p \lambda_2/\kappa \\ N \lambda_3 \\ \lambda_4
       \end{matrix}\right)\Bigr]^{-s}
  +
  \hspace{-0.5em}
  \sum_{\lambda \in \ZZ^4 \setminus \{0\}}
  \hspace{-0.5em}
  \chi(\lambda_4)
  P_Z \Bigl[\left(\begin{matrix}
       N p \lambda_1 \\ N^2 p \lambda_2/\kappa \\ N p \lambda_3 \\ \lambda_4
       \end{matrix}\right)\Bigr]^{-s}
\text{,}
\end{align*}
yielding \eqref{eq:firsteisensteinequality}.  The second equality \eqref{eq:secondeisensteinequality} follows from
\begin{align*}
&\quad
   \big( N^2 \pi^s \big)^s \Gamma(s)^{-1} \, \mathbb{E}^{(p) *} _{N p, N^2 p / \kappa, \chi} (Z, s)
\\[6pt]
&=
  p^{3s/2 }
  \hspace{-0.5em}
  \sum_{\lambda \in \ZZ^4 \setminus \{0\}}
  \hspace{-0.5em}
  \chi(\lambda_4) \,
  \Bigl(
  p^{-s}
  P_Z \Bigl[\left(\begin{matrix}
       \lambda_1/N  \\ p\lambda_2/\kappa  \\
       \lambda_3/N  \\ \lambda_4/N^2
       \end{matrix}\right)\Bigr]^{-s}
  +
  p^{-2 s}
  P_Z \Bigl[\left(\begin{matrix}
       \lambda_1/N \\ \lambda_2/\kappa  \\
       \lambda_3/N  \\ \lambda_4/N^2 p 
       \end{matrix}\right)\Bigr]^{-s}
  \Bigr)
\\[6pt]
&=
  p^{s/2} \hspace{-1em}
  \sum_{\substack{
        \alpha (N),\, \beta (\kappa) \\[1pt]
        \gamma (N),\, \delta (N^2)
       }}
  \hspace{-1em}
  \chi(\delta) \,
  \Bigl(
  p^{-s}
  \sum_{g (p)}
  \zeta\bigl(
        s,\, 0,\,
        (\frac{\alpha}{N},\, \frac{\beta}{\kappa},\,
         \frac{\gamma}{N},\, \frac{p \delta + N^2 g}{N^2 p})^\tr,\,
        P_Z
       \bigr)
\\[4pt]
& \hspace{9em}
  +
  p^{-1}
  \sum_{h (p)}
  \zeta\bigl(
        s,\, (0,\,\frac{\kappa h}{p},\, 0,\, 0)^\tr,\,
        (\frac{\alpha}{N},\, \frac{p \beta}{\kappa},\,
              \frac{\gamma}{N},\, \frac{\delta}{N^2})^\tr\hspace{-0.3em},\,
        P_Z
       \bigr)
  \Bigr)
\text{.}
\end{align*}
The last equality \eqref{eq:thirdeisensteinequality} can be deduced analogously.

The absolute convergence of \eqref{eisensteinseriesdef} and \eqref{eq:firsteisensteinequality}, follows from \eqref{eq:secondeisensteinequality} when reading the equality backwards.  The meromorphic continuation can be deduced from the same equality.  In the meromorphic continuation of the Epstein $\zeta$-function, at most a simple pole at $s = 2$ can occur (see~\citep{terrasharmonicanalysis}).  In the light of \eqref{eq:secondeisensteinequality}, the same holds for the meromorphic continuation of $\mathbb{E}^{(p) *} _{N p, N^2 p / \kappa, \chi} (Z, s)$.  The residue at $s = 2$ equals
\begin{gather*}
  N^{-4}
  \sum_{\substack{
        \alpha(N),\, \beta (\kappa) \\[1pt]
        \gamma (N),\, \delta (N^2)
       }}
  \hspace{-1em}
  \chi(\delta)
  (1 + 1)
=
  \begin{cases}
  2 \kappa \varphi(N)/N\text{,}  &\text{if $\chi$ is trivial;} \\
  0\text{,}                      &\text{otherwise.}
  \end{cases}
\end{gather*}
This completes the proof.
\end{proof}

\subsection{The functional equation}

We will establish the functional equation of $\mathbb{E}^{(p) *} _{N p, N^2 p / \kappa, \chi} (Z, s)$.  As a first step, we deduce a representation in terms of $P_Z$, which was defined in \eqref{eq:PZdefinition}.

Let $\chi_L$ be the primitive character with period $L$ that induces $\chi$.  Moreover, let $L R$ be the minimal period of $\chi$.  Recall that $\kappa \isdiv N$.
\begin{proposition}
\label{compeisensteinfcteqPZdecomp}
We have
\begin{gather}
  \mathbb{E}^{(p) *} _{N p, N^2 p / \kappa, \chi} (Z, 2 - s)
=
  \frac{\varphi(R)}{R}
  \frac{\kappa}{G_{\overline {\chi_L}}} \;
  \sum_{r \isdiv R}
   \big( (L r)^2 / \pi \big)^s \, \Gamma(s)  \,
   \chi_L (r) \,
   \frac{\mu(r)}{\varphi(r)}
   \sum_{\lambda \in \ZZ^4 \setminus \{0\}} \!\!
   ({\overline {\mathbbm{1}_r \chi_L}}) (\lambda_4) \,
   {\frak s}(r, \lambda)
\text{,}
\end{gather}
where
\begin{gather*}
  {\frak s}(r, \lambda)
:=
  \Bigl(
  p^{3s/2}
  P_{W_{N p} \langle Z\rangle}
  \Bigl[\left(\begin{matrix}
   L r p \lambda_1 \\ L r \kappa p \lambda_2 \\
   L r p \lambda_3 \\ \lambda_4
  \end{matrix}\right)\Bigr] ^{-s}
  +
  p^{s/2}
  P_{W_{N p} \langle Z\rangle}
  \Bigl[\left(\begin{matrix}
   L r \lambda_1 \\ L r \kappa p \lambda_2 \\
   L r \lambda_3 \\ \lambda_4
  \end{matrix}\right)\Bigr]^{-s}
  \Bigr)
\text{.}
\end{gather*}
\end{proposition}

If we impose an additional condition on $\kappa$, we can deduce a more convenient expression, based on the representation of the Eisenstein series given in Lemma~\ref{la:eisensteinreformulations}.
\begin{corollary}
\label{compeisensteinfcteq}
Suppose that $\kappa \isdiv L$.  Then $  \mathbb{E}^{(p) *} _{N p, N^2 p / \kappa, \chi} (Z, 2 - s)$ satisfies the functional equation
\begin{gather*}
  \mathbb{E}^{(p) *} _{N p, N^2 p / \kappa, \chi} (Z, 2 - s)
=
  \frac{\varphi(R)}{R}
  \frac{\kappa}{G_{\overline {\chi_L}}} \,
  \sum_{r \isdiv R}
   \chi_L (r)
   \frac{\mu(r)}{\varphi(r)} \,
   \mathbb{E}^{(p) *} _{L r p, L r \kappa p, \overline{\mathbbm{1}_r \chi_L}}
             (W_{N p} \langle Z\rangle, s)
\text{.}
\end{gather*}
\end{corollary}

\begin{proof}[{\tit Proof of Proposition \ref{compeisensteinfcteqPZdecomp}}]
Using \eqref{eq:secondeisensteinequality}, we find
\begin{align*}
& \quad
  \mathbb{E}^{(p) *} _{N p, N^2 p/ \kappa, \chi} (Z, 2 - s)
\\[6pt]
&=
  N^{-4} \big( N^2 / \pi^{-s} \big) \, \Gamma(s)
  \hspace{-1em}
  \sum_{\substack{
        \alpha (N),\, \beta (\kappa) \\[1pt]
        \gamma (N),\, \delta (N^2)
       }}
  \hspace{-1em}
  \chi(\delta)
\\[4pt]
& \hphantom{=}\,
  \cdot
  \Big(
  p^{-s/2}
  \hspace{-2em}
  \sum_{\substack{
        h (p) \\[1pt]
         v
        =
         \left(0,\,\kappa h/p,\, 0,\, 0\right)^\tr
       }}
  \;
  \sum_{\substack{
        \lambda \in \ZZ^4 \\[1pt]
        \lambda + v \ne 0
       }}
  \!
  e \big(
   -(\frac{\alpha}{N},\, \frac{p \beta}{\kappa},\,
     \frac{\gamma}{N},\, \frac{\delta}{N^2}) \,
   \lambda
  \big) \,
  P_Z ^{-1}
  [\lambda + v] ^{-s}
\\[4pt]
& \hspace{2.5em}
  +
  p^{-(2 - s)/2}
  \sum_{g (p)}\;
  \sum_{\lambda \in \ZZ^4 \setminus \{0\}}
  e \big(
   -(\frac{\alpha}{N},\, \frac{\beta}{\kappa},\,
     \frac{\gamma}{N},\, \frac{p \delta + N^2 g}{N^2 p}) \,
   \lambda
  \big) \,
  P_Z ^{-1}
  [\lambda] ^{-s}
  \Bigr)
\text{.}
\end{align*}
Set $\nu := N / L R$.  As a next step, evaluate the character sum and then replace $\lambda_1$ by $N \lambda_1$, $\lambda_2$ by $\kappa \lambda_2$, $\lambda_3$ by $N \lambda_3$, and $\lambda_4$ by $N \nu \lambda_4$.  In the next formula we write $r$ for $\gcd(R, \lambda_4)$:
\begin{align*}
&\quad
  \Bigl(
   \pi^{-s} \, \Gamma(s) \, \frac{\kappa}{R G_{{\overline \chi}_L}}
  \Bigr)^{-1} \,
  \mathbb{E}^{(p) *} _{N p, N^2 p / \kappa, \chi} (Z, 2 - s)
\\&
=
  N^{2 s} \chi_L (R)
  \Bigl(
  p^{- s/2}
  \hspace{-3em}
  \sum_{\substack{
        \lambda \in \ZZ^4 \\[1pt]
        h (p) \\[1pt]
        v = \left(N \lambda_1,\, \kappa \lambda_2,\, N \lambda_3,\,N \nu \lambda_4\right) \\[1pt]
        v + \left(0,\, \kappa h/p, 0, 0\right) \neq 0
       }}
  \hspace{-3em}
  {\overline {\chi_L}} (\lambda_4) \,
  \mu(r) \varphi(R / r) \,
  P_Z ^{-1}
  \Bigl[\left(\begin{matrix}
   N \lambda_1 \\ \kappa \lambda_2 + \kappa h/p \\
   N \lambda_3 \\ N \nu \lambda_4
  \end{matrix}\right)\Bigr] ^{-s}
\\[4pt]
 & \hspace{6.5em}
  +
  p^{s/2 - 1}
  \hspace{-1em}
  \sum_{\substack{
        \lambda \in \ZZ^4 \setminus \{0\} \\[1pt]
        g (p)
       }}
  \hspace{-0.5em}
  {\overline {\chi_L}} (\lambda_4) \,
  \mu(r) \varphi(R / r) \,
  e\big( -g N \nu \lambda_4 / p \big) \,
  P_Z ^{-1}
  \Bigl[\left(\begin{matrix}
   N \lambda_1 \\ \kappa \lambda_2 \\
   N \lambda_3 \\ N \nu \lambda_4
  \end{matrix}\right)\Bigr]^{-s}
  \Bigr)
\text{.}
\end{align*}
Fixing any $r$, we consider the partial sum
\begin{align*}
&\quad
  N^{2 s} \, \chi_L (R) \,
  \mu(r) \varphi(R / r)
\\[4pt] &
\hspace{4em}
  \cdot
  \Bigl(
  p^{- s/2}
  \hspace{-3em}
  \sum_{\substack{
        \lambda \in \ZZ^4 \\[1pt]
        h (p) \\[1pt]
        v = \left(N \lambda_1,\, \kappa \lambda_2,\, N \lambda_3,\,N \nu \lambda_4\right) \\[1pt]
        v + \left(0,\, \kappa h/p, 0, 0\right) \neq 0  \\[1pt]
        \gcd(R, \lambda_4) = R/r
       }}
  \hspace{-2.5em}
  {\overline {\chi_L}} (\lambda_4) \,
  P_{W_{N p} \langle Z\rangle}
  \Bigl[\left(\begin{matrix}
   N \lambda_1 \\ N \kappa p \lambda_2 + N \kappa h \\
   N \lambda_3 \\ \nu \lambda_4/p
  \end{matrix}\right)\Bigr] ^{-s}
\\ & \hspace{6.5em}
  +
  p^{s/2 - 1}
  \hspace{-1.5em}
  \sum_{\substack{
        \lambda \in \ZZ^4 \setminus \{0\} \\[1pt]
        g (p) \\[1pt]
        \gcd(R, \lambda_4) = R/r
       }}
  {\overline {\chi_L}} (\lambda_4) \,
  e\big(- g N \nu \lambda_4 / p \big) \,
  P_{W_{N p} \langle Z\rangle}
  \Bigl[\left(\begin{matrix}
   N \lambda_1 \\ N \kappa p \lambda_2 \\
   N \lambda_3 \\ \nu \lambda_4/p
  \end{matrix}\right)\Bigr]^{-s}
  \Bigr)
\text{.}
\end{align*}
We evaluate the sum over $g$, substitute $p R\lambda_4/r$, in the first sum, and $R \lambda_4/r$, in the second sum, by $\lambda_4$, and replace $p \lambda_2 + h$ by $\lambda_2$.  Recalling that $p \equiv 1\;(N)$, we find that the above expression equals
\begin{align*}
&
  (L r)^{2 s} \chi_L (r) \,
  \mu(r) \varphi(R / r) \,
  \sum_{\lambda \in \ZZ^4 \setminus \{0\}}
  ({\overline {\mathbbm{1}_r \chi_L}}) (\lambda_4) \cdot
\\ & \hspace{7em}
  \cdot
  \Bigl(
  p^{3s/2}
  P_{W_{N p} \langle Z\rangle}
  \Bigl[\left(\begin{matrix}
   L r p \lambda_1 \\ L r \kappa p \lambda_2 \\
   L r p \lambda_3 \\ \lambda_4
  \end{matrix}\right)\Bigr] ^{-s}
  +
  p^{s/2}
  P_{W_{N p} \langle Z\rangle}
  \Bigl[\left(\begin{matrix}
   L r \lambda_1 \\ L r \kappa p \lambda_2 \\
   L r \lambda_3 \\ \lambda_4
  \end{matrix}\right)\Bigr]^{-s}
  \Bigr)
\text{.}
\end{align*}
Since $R$ is square free, we have $\varphi\left(R/r\right) = \varphi(R)/\varphi(r)$.  Summing the last expression for all $r \isdiv R$ yields the claim.
\end{proof}

In general, the Eisenstein series occurring on the right hand side of the functional equation in Corollary \ref{compeisensteinfcteq} are hard to treat. After summing over divisors of $R$, we can reduce them to a single one.

\begin{proposition}
\label{compeisensteinseriesfcteqrrhs}
If $\kappa \isdiv L$, then
\begin{gather}
\label{eq:compeisensteinseriesfcteqrrhs}
  \sum_{\theta \isdiv R}
  \frac{\mu (\theta) \varphi (\theta)}{\theta} \,
  \mathbb{E}^{(p) *} _{L R p / \theta, L^2 R^2 p / (\kappa \theta^2),
                       \overline{\mathbbm{1}_{R/r} \chi_L}}
            (W_{p L R/\theta} \langle Z\rangle, 2 - s)
=
  \chi_L (R)
  \frac{\mu (R)}{R} \,
  \frac{\kappa}{G_{\overline{\chi_L}}} \,
  \mathbb{E}^{(p) *}_{L R p, L R \kappa p, \overline{\mathbbm{1}_R \chi_L}}
  \left(Z,\,s\right)
\text{.}
\end{gather}
\end{proposition}

\begin{proof}
Use the functional equation in Corollary \ref{compeisensteinfcteq} and the equality
\begin{multline*}
  \sum_{\theta \isdiv R}
  \frac{\mu (\theta) \varphi (\theta)}{\theta}
  \frac{\varphi\left(R/\theta\right)}{R/\theta}
  \frac{\kappa}{G_{\overline{\chi_L}}}
  \sum_{r \isdiv \gcd(R/\theta)} \!
  \chi_L (r)
  \frac{\mu(r)}{\varphi(r)} \,
  \mathbb{E}^{(p) *} _{L r p,\,L \kappa r p,\,\overline{\mathbbm{1}_r \chi_L}}
  \left(Z,\,s\right)
\\[4pt]
=
  \frac{\varphi(R)}{R}
  \frac{\kappa}{G_{\overline{\chi_L}}}
  \sum_{r \isdiv R} \,
  \Bigl(
   \sum_{\theta \isdiv \gcd(R/r)}
   \mu(\theta)
  \Bigr)\,
  \chi_L(r)
  \frac{\mu(r)}{\varphi(r)} \,
  \mathbb{E}^{(p) *} _{L r p,\,L \kappa r p,\,\overline{\mathbbm{1}_r \chi_L}}
  \left(Z,\,s\right)
\text{.}
\qedhere
\end{multline*}
\end{proof}

In Section \ref{sec:fstintegralrepresentation}, we will deal with an integral involving the left hand side of \eqref{eq:compeisensteinseriesfcteqrrhs}.  Lemma \ref{la:eisensteinreductionproperties} provides basic tools to investigate it based on mathematical induction.  Let $q$ be a prime and $r \in \NN$ such that $\gcd(p, N q r) = 1$.  Given $Z \in \HS_2$ and $s\in \CC$ satisfying $\Re(s) > 2$, define
\begin{gather}
\label{eq:eisensteindifference_definition}
  \mathcal{E}^{(p) *}_{N p, N^2 p / \kappa, q, r, \chi}
  (\cdot,\,s)
  \left| W_{N q r p} \right.
:=
  \mathbb{E}^{(p) *}_{N p,\, N^2 p / \kappa,\, \chi}
  (\cdot,\,s)\left|W_{N p}\right.
  -
  \mathbb{E}^{(p) *}_{N q p,\, (N q)^2 p / \kappa,\, \mathbbm{1}_q \chi}
  (\cdot,\,s)\left|W_{N q p}\right.
\text{.}
\end{gather}

\begin{lemma}
\label{la:eisensteinreductionproperties}
Given $Z \in \HS_2$ and $s\in \CC$ satisfying $\Re(s) > 2$, the following statements are true.
\begin{itemize}
\item[(i)] The considered difference of Eisenstein series satisfies the equality
      \begin{multline*}
        \mathcal{E}^{(p) *}_{N p, N^2 p / \kappa, q, r, \chi} (Z,\,s)
      =
        \big( (N r q)^2 / \pi \big)^s \, \Gamma(s)
        \sum_{\lambda \in \ZZ^4 \setminus \{0\}}
        \chi (q \lambda_4) \,\cdot
        \hspace{3em}
      \\[4pt]
        \cdot
        \Bigl(
        p^{s / 2}
        P_Z
        \Bigl[\left(\begin{matrix} 
         N q r \lambda_1 \\
         N^2 q^2 r^2 p \lambda_2 / \kappa \\
         N q r \lambda_3 \\
         q \lambda_4
        \end{matrix}\right)\Bigr]^{-s}
        +
        p^{3 s / 2}
        P_Z
        \Bigl[\left(\begin{matrix} 
         N q r p \lambda_1 \\
         N^2 q^2 r^2 p \lambda_2 / \kappa \\
         N q r p \lambda_3 \\
         q \lambda_4
        \end{matrix}\right)\Bigr]^{-s}
        \Bigr)
      \text{.}
     \end{multline*}\\[2pt]
\item[(ii)] \label{it:eisensteindifferencesecondproperty} If $q \isdiv N$, we have $\mathcal{E}^{(p) *}_{N p, N^2 p / \kappa, q, r, \chi} (Z,\,s) = 0$. \\[2pt]
\item[(iii)] The function $\mathcal{E}^{(p) *}_{N p, N^2 p / \kappa, q, r, \chi} (\cdot,\,s)$ is a Siegel modular form of weight $0$ with character $\ov{\chi}$ for $\Gamma_{2,1} \big( N q r p, (N q r p)^2 \big)$. \\[2pt]
\item[(iv)] We have $
        \mathcal{E}^{(p) *}_{N p, N^2 p / \kappa, q, r, \chi} (\cdot,\,s)
        \left| M_{N^2 q r^2 p / \kappa} ^\tr \right.
      = \mathcal{E}^{(p) *}_{N p, N^2 p / \kappa, q, r, \chi} (\cdot,\,s)
      $.
\end{itemize}
\end{lemma}

\begin{proof}
In order to prove (i), we consider
\begin{align*}
&\quad
  \Big( \big(N^2 / \pi \big)^s \, \Gamma(s) \Big)^{-1}
  \mathbb{E}^{(p) *} _{N p,\, N^2 p / \kappa,\, \chi}
  \left(
   W_{N p}\langle Z\rangle,\, s
  \right)
\\[6pt]
&=
  \sum_{\lambda \in \ZZ^4 \setminus \{0\}}
  \hspace{-0.5em}
  \chi(\lambda_4)
  \Bigl(
   p^{s/2}
   P_{W_{N p}\langle Z\rangle}
   \Bigl[\left(\begin{matrix}
    N \lambda_1 \\
    N^2 p \lambda_2/ \kappa \\
    N \lambda_3 \\
    \lambda_4
   \end{matrix}\right)\Bigr]^{-s}
   +
   p^{3s/2}
   P_{W_{N p}\langle Z\rangle}
   \Bigl[\left(\begin{matrix}
    N p \lambda_1 \\
     N^2 p \lambda_2/ \kappa \\
    N p \lambda_3 \\
    \lambda_4
   \end{matrix}\right)\Bigr]^{-s}
  \Bigr)
\\[6pt]
&=
  q^{2 s}
  \hspace{-0.5em}
  \sum_{\lambda \in \ZZ^4 \setminus \{0\}}
  \hspace{-0.5em}
  \chi(\lambda_4) \cdot
\\[4pt]
& \hspace{3em}
  \cdot
   \Bigl(
   p^{s/2}
   P_{W_{N q p}\langle Z\rangle}
   \Bigl[\left(\begin{matrix}
    N q \lambda_1 \\
     N^2 q^2p \lambda_2/ \kappa \\
    N q \lambda_3 \\
    \lambda_4
   \end{matrix}\right)\Bigr]^{-s}
   +
   p^{3 s/2}
   P_{W_{N q p}\langle Z\rangle}
   \Bigl[\left(\begin{matrix}
    N q p\lambda_1 \\
     N^2 q^2p \lambda_2/ \kappa \\
    N q p \lambda_3 \\
    \lambda_4
   \end{matrix}\right)\Bigr]^{-s}
  \Bigr)
\text{.}
\end{align*}
We split the addends with respect to $q \isdiv \lambda_4$ and $q \nisdiv \lambda_4$.
\begin{align*}
& \quad
  \big(\pi^{-s} \, \Gamma(s) \big)^{-1}
  \mathbb{E}^{(p) *} _{N p,\, N^2 p / \kappa,\, \chi}
  \left(
   W_{N p}\langle Z\rangle,\, s
  \right)
\\[6pt]
&=
  \big( \pi^{-s} \, \Gamma(s) \big)^{-1}
  \mathbb{E}^{(p) *} _{N q p,\, N^2 q^2 p / \kappa,\, \mathbbm{1}_q \chi}
  \left(
   W_{N q p}\langle Z\rangle,\, s
  \right)
\\[4pt]
&\qquad
  +
  (N q)^{2 s}
  \sum_{\lambda \in \ZZ^4 \setminus \{0\}}
  \chi(q \lambda_4) \cdot
\\[4pt]
& \hspace{5em}
  \cdot
  \Bigl(
   p^{s/2}
   P_{W_{N q p}\langle Z\rangle}
   \Bigl[\left(\begin{matrix}
    N q \lambda_1 \\
     N^2 q^2 p \lambda_2/ \kappa \\
    N q \lambda_3 \\
    q \lambda_4
   \end{matrix}\right)\Bigr]^{-s}
   +
   p^{3s/2 }
   P_{W_{N q p}\langle Z\rangle}
   \Bigl[\left(\begin{matrix}
    N q p \lambda_1 \\
     N^2 q^2 p \lambda_2/ \kappa \\
    N q p \lambda_3 \\
    q \lambda_4
   \end{matrix}\right)\Bigr]^{-s}
  \Bigr)
\text{.}
\end{align*}

This yields (i).  As an immediate consequence, we deduce (ii), since $\chi(q \lambda_4) = 0$ as $1 \ne q \isdiv N$.  Using the defining equation \eqref{eq:eisensteindifference_definition} of $\mathcal{E}^{(p) *}$, we prove (iii).  The last claim (iv) is a consequence of (i).
\end{proof}


\section{A first integral representation for the Rankin convolution}
\label{sec:fstintegralrepresentation}

Recall the notation $e(x) = e^{2 \pi i x}$, which we will use frequently in the next two sections.  Also recall the notation for the elements of $\HS_2$ given in~\eqref{eq:HSelements_notation}.  The Petersson scalar product for Jacobi forms $\phi$ and $\psi$ of weight $k$ and index $m$ is
\begin{gather}
  \langle \phi,\, \psi \rangle
:=
  \int_{\cF^*} \phi(\tau, z)\, \psi(\tau, z)\,
                 y^{k-3} \, e^{-4 \pi m y^{-1} v^2} \; dx\, dy\, du\, dv
\text{,}
\end{gather}
where $\cF^*$ is a fundamental domain with respect to the action of the full Jacobi group $\SL{2}(\ZZ) \ltimes \ZZ^2$ on $\HS_1 \times \CC$.  The reader is referred to \citep{EZ85} for details.

A fundamental domain for the action of $\Gamma^{(p)}_{2,1}$ on $\HS_2$ can be constructed from this as follows:
\begin{gather}
  \cF^{(p)}_{2,1}
:=
  \Big\{
   Z \in \HS_2 \,:\, (\tau, z) \in \cF^*,\, x' \in [0, 1/p],\, y' > y^{-1} v^2
  \Big\}
\text{.}
\end{gather}

Let $f$ and $g$ paramodular forms in $\big[\Gamma^{(p)}_{2,1}(p, p), k, 1 \big]_0$ with Fourier-Jacobi expansions 
\begin{gather}
\label{eq:fourierjacobicoefficientsdefinition}
  f(Z)
=
  \sum_{m=1} ^\infty f_m(\tau, z)\, e^{2 \pi i\, m \tau'}
\quad\text{and}\quad
  g(Z)
=
  \sum_{m=1} ^\infty g_m(\tau, z)\, e^{2 \pi i\, m \tau'}
\text{.}
\end{gather}
We attach a Dirichlet series, the Rankin convolution, to $f$ and $g$:
\begin{gather}
  D_{f,g} (s)
:=
  \sum_{m = 1} ^\infty \langle f_m, g_m \rangle \, m^{-s}
\text{,}
\end{gather}
for all $s$ such that the right hand side converges locally uniformly absolutely.  The twisted Rankin convolution and its completion are
\begin{align}
  D_{f,g,\chi} (s)
&:=
  \sum_{m = 1} ^\infty \chi(m) \langle f_{m}, g_{m} \rangle m^{-s},
\\[6pt]
  \mathbb{D}_{f,g,\chi} (s)
&:=
  ( 2 \pi/N )^{-2 s} \,
  \Gamma(s) \Gamma(s - k + n) \,
  L(2 s - 2 k + 2 n, \chi^2) \,
  D_{f, g, \chi} (s)
\text{.}
\end{align}

We will write
\begin{gather}
  \big\{f,\,g \big\}
:=
  \frac{1}{\big[\Gamma_{2,1}^{(p)*}(p,p) \,:\, \Gamma'\big]}
  \int_{\cF'} f(Z)\, \ov{g(Z)}\, \det(Y)^{k-3}
             \; \,dX\,dY
\end{gather}
for the Petersson scalar product of $f$ and $g$.  Here, $\Gamma' \subseteq \Gamma_{2,1}^{(p)*}(p,p)$ is any group that $f$ and $g$ are invariant under, and $\cF'$ is an associated fundamental domain.

\begin{theorem}
\label{fstintegralrepresentation}
Let $\chi^{(2)}$ be a primitive character that induces $\chi^2$.  Let $L$ denote the minimal period of $\chi^{(2)}$.  Fix $\theta \in \NN$ satisfying $\theta \isdiv N/L^{(2)}$.  If $\theta \ne 1$ assume, in addition, that $\chi$ is primitive. Then we have
\begin{multline}
\label{eq:fstintegralrepresentation}
  \Big\{
   \mathbb{E}^{(p) *} _{N p / \theta,\, N^2 p / \theta^2,\,
                        \mathbbm{1}_{N / \theta} \chi^{(2)}}
   (\cdot,s)
   \big| W_{N p / \theta} \,
   f_\chi \big| W_{N p} ,\,
   g\big|W_{N p}
  \Big\}
\\[6pt]
=
  \frac{2}{i_N}
  ( \pi / N^2 )^{k-2} \,
  p^{3s/2 - 1}
  (1 + p^{-s}) \,
  \mathbb{D}_{f,\,g,\,\chi} (s + k - 2)
\text{.}
\end{multline}
\end{theorem}

\begin{proof}
The first and second argument of the Petersson scalar product on the left hand side of \eqref{eq:fstintegralrepresentation} have the common symmetry group $\Gamma_{2,1} \big(N p, (N p)^2 \big)$.  Since $W_{N p}$ normalizes this group, we see that the scalar product converges absolutely due to the polynomial growth of the Eisenstein series and
\begin{gather*}
  \Big\{
   \mathbb{E}^{(p) *} _{ N p / \theta,\, N^2 p / \theta^2 ,\,
                        \mathbbm{1}_{N/\theta} \chi^{(2)}}
   (\cdot,s)
   \big| W_{N p / \theta} \,
   f_\chi \big| W_{N p} ,\,
   g \big|W_{N p}
  \Big\}
=
   \Big\{
   \mathbb{E}^{(p) *} _{N p / \theta,\, N^2 p / \theta^2,\,
                        \mathbbm{1}_{N/\theta} \chi^{(2)}}
    (\cdot,s) \big| D_{1/\theta} \,
   f_\chi, \,
   g
  \Big\}
\text{.}
\end{gather*}

We use mathematical induction on the prime divisors of $\theta$.  If $\theta = 1$, we can apply the unfolding trick to \eqref{eq:fstintegralrepresentation}.  The above scalar product simplifies as follows:
\begin{align*}
& \quad
  \Big(
  \frac{2}{i_N}
  \big( N^2 / \pi \big)^s \,
  p^{3s/2}
  (1 + p^{-s}) \,
  \Gamma(s)
  L(2 s, \chi ^2)
  \Big)^{-1}
  \Big\{
   \mathbb{E}^{(p) *} _{N p,\, N^2 p,\,
                        \mathbbm{1}_{N} \chi^{(2)}}
   (\cdot,s) \,\
   f_\chi ,\,
   g
  \Big\}
\\[6pt]
&=
  \int_{\cF^{(p)} _{2,1}}
  f_\chi (Z) \, \ov{g(Z)} \, y^{k - 3}
  \big(y' - y^{-1} v^2 \big)^{s + k - 3}
  \; dX\, dY
\displaybreak[0]
\\[6pt]
&=
  \int_{\mathcal{F}^{(p)}}
  \;
  \sum_{m_f,\,m_g = 1} ^{\infty}
  \chi(m_f) \,
  f_{p m_f} (\tau, z) \,
  \overline{g_{p m_g}(\tau, z)} \,
  y^{k - 3}
  \big( y' - y^{-1} v^2 \big)^{s + k - 3}
\\[4pt]
&\hspace{12em}
  \cdot
   e\big( i p (m_f + m_g) y' + p (m_f - m_g) x' \big)
  \; dX\, dY
\\[6pt]
&=
  p^{-1}
  (4 \pi)^{-s - k + 2}
  \Gamma(s + k - 2)
\\[4pt]
&\qquad
  \cdot
  \int_{(Z_1,\,w) \in \cF ^*}
  \;
  \sum_{m = 1} ^{\infty}
  \chi(m)
  (p m)^{-s - k + 2} \,
  f_{p m} (\tau, z) \,
  \overline{g_{p m}(\tau, z)} \,
  e^{- 4 \pi p m y^{-1} v^2}
\\[4pt]
&\hspace{12em}
  \cdot
  y^{k - 3}
  \; du\, dv\, dx\, dy
\\[6pt]
&=
  (4 \pi)^{-s - k + 2} \,
  p^{-1} \,
  \Gamma(s + k - 2) \,
  D_{f,\,g,\,\chi} (s + k - 2)
\text{.}
\end{align*}
This proves the statement if $\theta = 1$.

Assume that $\theta \ne 1$, and let $q \isdiv \theta$ be a prime. Recall the definition of $\mathcal{E}^{(p)*}$ in \eqref{eq:eisensteindifference_definition}.  It suffices to show that
\begin{align*}
&\quad
  \Big\{
   \mathcal{E}^{(p) *}_{N p / \theta, N^2 p / \theta^2, q, \theta/q, \chi^2}
   (\cdot,\,s)
   \big| W_{N p} \,
   f_\chi \big| W_{N p} ,\,
   g\big|W_{N p}
  \Big\}
=
  0
\text{.}
\end{align*}
By (ii) of Proposition \ref{compeisensteinseriesfcteqrrhs}, this is true, whenever $q \isdiv N/\theta$.  Hence we may assume that $q \nisdiv N / \theta$.

The modular forms on the left hand side of the last equation share the symmetry group $\Gamma_{2,1} ^1 \big(N p,\, (N p)^2\big)$, and $M_{N p} ^\tr$ normalizes this group.  Hence so does $M_{(N p)^2 \nu/q} ^\tr$ for any $\nu \in \NN$.  Using the invariance of $g$ and $\mathcal{E}^{(p)*}$ under $M_{(N p)^2 \nu/q} ^\tr$, we find
\begin{gather*}
  \Big\{
   \mathcal{E}^{(p) *}_{N p / \theta, N^2 p / \theta^2, q, \theta/q, \chi^2}
   (\cdot,\,s) \,
   f_\chi ,\,
   g
  \Big\}
=
  q^{-1} \,
  \sum_{\nu (q)}
  \Big\{
   f_\chi \big| M_{(N p)^2 \nu/q}^\tr \,
   \mathcal{E}^{(p) *}_{N p / \theta, N^2 p / \theta^2, q, \theta/q, \chi^2}
   (\cdot,\,s) ,\,
   g
  \Big\}
\text{.}
\end{gather*}
This sum vanishes by Lemma \ref{la:vanishingundertranslations}, yielding the claim.
\end{proof}

An immediate consequence of Theorem \ref{fstintegralrepresentation} is
\begin{corollary}
\label{fstintegralrepresentationsum}
Suppose that $\chi$ is primitive. Let $\chi^{(2)}$ be a primitive character that induces $\chi^2$.  Denote the minimal period of $\chi^{(2)}$ by $L$.  Let $L R$ be the minimal period of $\chi^2$. Then we have
\begin{multline*}
  \sum_{\theta \mid R}   
  \frac{\mu(\theta) \varphi(\theta)}{\theta}
  \Big\{
  \mathbb{E}^{(p) *}_{L R p / \theta,\,
                       \left( L R p / \theta \right)^2,\,
                       \mathbbm{1}_{N/\theta} \chi^{(2)}}
   (\cdot,s)
   \big| W_{L R / \theta} \,
   f_{\chi} \big| W_{N p} ,\,
   g \big| W_{N p}
  \Big\}
\\[6pt]=
  \frac{2}{i_N R^{(2)}}
  \big( \pi / N^2 )^{k - 2} \,
  p^{3s/2 - 1}
  (1 + p^{-s}) \,
  \mathbb{D}_{f,\,g,\,\chi} (s + k - 2)
\text{.}
\end{multline*}
\end{corollary}


\section{A second integral representation for the Rankin convolution}
\label{sec:sndintegralrepresentation}

In this section, we assume that $f, g \in \big[\Gamma^{(p)}_{2,1}(p, p), k, 1\big]_0$.  In particular, the Fourier-Jacobi coefficients $f_m$ and $g_m$ of $f$ and $g$, defined in \eqref{eq:fourierjacobicoefficientsdefinition}, vanish, if $p \nisdiv m$.  Throughout the section, fix $\chi^{(2)}$, $L$ and $R$ as in Corollary \ref{fstintegralrepresentationsum}.  In addition, set $\nu := N / L R$.  Note that $\nu \isdiv 4$.

\begin{lemma}
\label{sndintreprfourierjacobiseries}
Suppose that $d \isdiv \nu$ and $C \isdiv N \nu$.  Moreover, let $\mu$ be the minimal, positive integer dividing $\nu$, such that $C \isdiv N d \mu$.  Choose any $\gamma, \gamma^*, p^*, \mu^* \in \NN$ such that $\gcd(\gamma, N \nu / C) = 1$, $\gamma \gamma^* \equiv 1 \; (N \nu / C)$, $p p^* \equiv 1 \; (N \nu / C)$ and $\mu \mu^* \equiv 1 \; (N d \mu p / C)$.  Fix $\epsilon \in \{\pm 1 \}$ and $\lambda \in \ZZ^\times$.  Set $\theta = N/\nu$, and choose $M_\lambda$ and $M_{d, C \gamma}$ as in Section~\ref{sec:paramodulargroup}.

For any $\epsilon_\gamma \equiv \gamma\; (\nu)$, we have the following Fourier-Jacobi expansion:
\begin{gather*}
  \Big(f \big| W_{N \nu p} M_{d, C \gamma} M_\lambda D_{\epsilon}\Big) (Z)
=
  \sum_{m = 1} ^\infty
  f_{m, C,d,\lambda, \epsilon, \epsilon_\gamma} (\tau, z) \,
  e\big( C^2  m \tau' / \mu^2 + C {\mu ^*}^2 \gamma^* p^* m / N \nu \big)
\text{.}
\end{gather*}
The function $f_{m, C,d,\lambda,\epsilon,\epsilon_\gamma}$ depends on its arguments and indices only.  If $p \nisdiv m$, we have $f_{m, C,d,\lambda,\epsilon,\epsilon_\gamma} = 0$.  If $\mu = \mu^* = 1$ and $C \nu \isdiv N$, we have
\begin{multline*}
  f_{m, C,d,\lambda, \epsilon, \epsilon_\gamma} (\tau, z)
=
  C^k
  f_m \Big(
       \tau,
       C z
       + p p^* \epsilon \epsilon_\gamma \big( \tau d \lambda_1 + d \lambda_2\big) / \nu
      \Big)
\\
  \cdot
  e\big( m (p p^* \epsilon \epsilon_\gamma/\nu )^2
     (\tau d^2 \lambda_1 ^2 + d^2 \lambda_1 \lambda_2) \big) \,
  e\big( 2 m C p p^* \epsilon \epsilon_\gamma w d \lambda_1 / \nu \big)
\text{.}
\end{multline*}
\end{lemma}
\begin{proof}
Set 
\begin{align*}
  H_1
&
:= \left(\begin{array}{cc|cc}
    1 &                                             &   &                 \\[2pt]
      & \frac{(1 - p \gamma p^* \gamma^*) C}{N \nu} &   & p^* \gamma^*    \\[2pt]
   \hline
      &                                             & 1 &                 \\[2pt]
      & - p \gamma                                  &   & \frac{N \nu}{C}
  \end{array}\right)
\text{,}
\\[6pt]
  H_2
&
:= \left(\begin{array}{cc|cc}
    \mu                                 & \frac{N d p \mu}{C} &  & \\[2pt]
    \frac{(\mu \mu^* - 1) C}{N d p \mu} & \mu^*               &  & \\[2pt]
   \hline
      &  & \mu^*                 & \frac{(1 - \mu \mu^*) C}{N d p \mu} \\[2pt]
      &  & - \frac{N d p \mu}{C} & \mu
   \end{array}\right)
\text{,}
\\[6pt]
   H_3
&
:= \left(\begin{array}{cc|cc}
    1 &   & \frac{N d^2 p^2 p^* \mu^2}{C} r & d p p^* \mu^* \mu r  \\[2pt]
      & 1 & d p p^* \mu^* \mu r             &                     \\[2pt]
    \hline
      &   & 1                               &                     \\[2pt]
      &   &                                 & 1
   \end{array}\right)
\text{,}
\quad
   H_4
:=
  \left(\begin{array}{cc|cc}
    1 &   & \frac{N d^2 p^2 p^* \epsilon_\gamma}{C \nu} &   \\
      & 1 &                                             &   \\
   \hline
       &   & 1                                          &   \\
       &   &                                            & 1
   \end{array}\right)
\text{,}
\end{align*}
where $r = (\epsilon_\gamma - \gamma^*) / \nu$.  The transformations $H_1$, $H_2$, and $H_3$ are elements of $\Gamma_{2, 1}(p, p)$.Recall that $W_p \in \Gamma^{(p)}_{2,1}(p,p)$.  Consequently,
\begin{gather*}
  f \big| W_{N \nu p} \, M_{d, C \gamma} \, M_\lambda \, D_{\epsilon}
=
  f \big| H_3 \, H_2 \, H_1 \, W_p \, W_{N \nu p} \, M_{d,C \gamma} \, M_\lambda \, D_\epsilon
\text{.}
\end{gather*}
Up to sign, the second transformation equals
\begin{gather*}
  \left(\begin{matrix}
  \mu                                 &               & \frac{N d^2 p^2 p^* \mu \, \epsilon_\gamma}{C \nu}  & \frac{d \, p \, p^* \mu \, \epsilon_\gamma}{C \nu} \\[8pt]
  \frac{(\mu \mu^* - 1) C}{N \, d \, p \, \mu} & \frac{C}{\mu} & \frac{d \, p \, p^* \mu^* \epsilon_\gamma}{\nu}        & \frac{p^* \mu^* (\epsilon_\gamma -  \nu \, r \, \mu \, \mu^*)}{N \nu} \\[8pt]
                                      &               & \frac{1}{\mu}                                  & \frac{1 - \mu \mu^*}{N \, d \, p \, \mu} \\[8pt]
                                      &               &                                                & \frac{\mu}{C}
  \end{matrix}\right)
\text{,}
\end{gather*}
yielding the first claim.

To prove the second claim assume that $\mu = \mu^* = 1$ and $C \nu \isdiv N$.  In this case we have \mbox{$H_4 \in \Gamma_{2,1} (p, p)$}, and thus
\begin{gather*}
  f \big| W_{N \nu p} \, M_{d, C \gamma} \, M_\lambda \, D_{\epsilon}
=
  f \big| D_\epsilon ^{-1} \, M_\lambda^{-1} \, H_4 \, H_3 \, H_2 \, H_1 \, W_p \, W_{N \nu p} \, M_{d,C \gamma} \, M_\lambda \, D_\epsilon
\text{.}
\end{gather*}
The second transformation equals
\begin{gather*}
  \left(\begin{matrix}
  1                                                             &   &                                                                & \frac{d \, p \, p^* \epsilon \, \epsilon_\gamma}{C \nu} \lambda_2 \\[8pt]
  \frac{d \, p \, p^* \epsilon \, \epsilon_\gamma}{\nu} \lambda_1 & C & \frac{d \, p \, p^* \epsilon \, \epsilon_\gamma}{\nu} \lambda_2 & \frac{p^* \gamma^*}{N \nu} \\[8pt]
                                                                &   & 1                                                              & - \frac{d \, p \, p^* \epsilon \, \epsilon_\gamma}{C \nu} \lambda_1 \\[8pt]
                                                                &   &                                                                & \frac{1}{C}
  \end{matrix}\right)
\text{.}
\end{gather*}
This proves the lemma.
\end{proof}

To prove Theorem \ref{sndintregralrepresentation}, we will need to evaluate a certain character sum.  The next Lemma was given in \citesndintkuss, and we reproduce the calculation for convenience only.  We will denote a sum ranging over representatives modulo $N$ that are units in $\ZZ / N \ZZ$ by $\beta \; (N)^\times$.  Observe that the quantity $A_{\chi, {\tilde \nu}}$ in \citep{kussthesis} features an additional factor $\nu^{-2}$.
\begin{lemma}
\label{Achisum}
Suppose that $\td \nu \isdiv \nu$.  Set $r := N/\gcd(N, L \nu)$.  Given $\beta,\, \gamma \in \ZZ$, choose $\beta^*,\, \gamma^* \in \ZZ$ satisfying \mbox{$\beta \beta^* \equiv \gamma \gamma^* \equiv 1\; (N \nu)$}.  For any $m \in \ZZ$, we have
\begin{align*}
  A_{\chi, {\tilde \nu}} (m)
&:=
\sum_{\substack{
      \beta (N {\tilde \nu}) ^{\times} \\
      \gamma (N {\tilde \nu}) ^{\times} \\
      \beta \equiv \gamma\; ({\tilde \nu})
     }}
  \chi
  \big(
  (\beta - \gamma) / {\tilde \nu}
  \big) \,
  e\big( m (\gamma^* - \beta^*)/N {\tilde \nu}\big)
\\[6pt]
&=
  \left\{
  \begin{aligned}
   &\chi(- m) \chi^{(2)}(R) \, \overline {\chi^{(2)}}\big(L R {\tilde \nu} / N \big) \,
   \mu(r) \varphi\big(R / r \big) \,
   \frac{\tilde \nu}{L R}  
   G_{\overline \chi} ^3 G_{\chi^{(2)}}
  \text{,}
             && \text{if } {\tilde \nu} = \nu \text{;} \\[4pt]
   &0
  \text{,}   && \text{otherwise.}
  \end{aligned}
  \right.
\end{align*}
In particular, we have $A_{\chi, {\tilde \nu}} (m) = \chi (m) \, A_{\chi, {\tilde \nu}} (1)$.
\end{lemma}

\begin{proof}
As a first step, we isolate a Gau\ss\ sum:
\begin{align*}
A_{\chi, \nu} (m)
&=
  \sum_{\substack{
      \beta (N {\tilde \nu}) ^{\times} \\
      \gamma (N {\tilde \nu}) ^{\times} \\
      \beta \equiv \gamma\; ({\tilde \nu})
     }}
  \chi(\beta \gamma)
  \chi \big( (\gamma^* - \beta^*) / {\tilde \nu} \big) \,
  e\big( m (\gamma^* - \beta^*) / N {\tilde \nu} \big)
\displaybreak[0]
\\[6pt]
&=
  \sum_{\substack{
      \beta (N {\tilde \nu}) \\
      \gamma (N {\tilde \nu}) \\
      \beta \equiv \gamma\; ({\tilde \nu})
     }}
  \overline \chi(\beta \gamma)
  \chi \big( (\gamma - \beta) / {\tilde \nu} \big) \,
  e\big( m (\gamma - \beta) / N {\tilde \nu} \big)
\displaybreak[0]
\\[6pt]
&=
  \sum_{\substack{
      \beta (N {\tilde \nu}) \\
      r (N)
     }}
  \overline \chi \big( \beta (\beta + {\tilde \nu} r) \big)
  \chi (r) \,
  e\big( m r/N \big)
\displaybreak[0]
\\[6pt]
&=
  {\tilde \nu}
  \sum_{\substack{
      \beta (N)^\times \\
      r (N)^\times
     }}
  \overline \chi \big( \beta (\beta + {\tilde \nu} r) \big)
  \chi (r) \,
  e\big( m r/N \big)
\displaybreak[0]
\\[6pt]
&=
  {\tilde \nu}
  \sum_{\substack{
      \beta (N)^\times \\
      r (N)^\times
     }}
  \overline \chi \big( \beta (\beta + {\tilde \nu}) \big)
  \overline \chi (r) \,
  e\big( m r/N \big)
\displaybreak[0]
\\[6pt]
&=
  {\tilde \nu}
  \chi(m) G_{\overline \chi}
  \sum_{\beta (N)^\times}
  \overline \chi \big( \beta (\beta + {\tilde \nu}) \big)
\text{.}
\end{align*}
Using $
  G_\chi \, \overline \chi (\beta + {\tilde \nu})
=
  \sum_{\gamma (N)^\times}
  \chi (\gamma) \, e\big( (\beta + {\tilde \nu})\, \gamma/N \big)
$, we obtain
\begin{align*}
  A_{\chi, \nu} (m)
&=
  \chi(m) \, {\tilde \nu} G_{\overline \chi}
  \sum_{\beta (N)^\times}
  \overline \chi \big( \beta (\beta + {\tilde \nu}) \big)
\displaybreak[0]
\\[6pt]
&=
  \chi(m) \, {\tilde \nu} \big( G_{\overline \chi} / G_\chi \big)
  \sum_{\substack{
        \beta (N)^\times \\
        \gamma (N)^\times
       }}
  \overline \chi(\beta)
  \chi (\gamma) \,
   e\big( (\beta + {\tilde \nu}) \gamma/N \big)
\displaybreak[0]
\\[6pt]
&=
  \chi(-m) \, {\tilde \nu} \big( G_{\overline \chi} ^2 / N \big)
  \sum_{\substack{
        \beta (N)^\times \\
        \gamma (N)^\times
       }}
  \overline \chi(\beta)
  \chi (\gamma) \,
  e\big( (\beta + {\tilde \nu}) \gamma/N \big)
\displaybreak[0]
\\[6pt]
&=
  \chi(-m) \, {\tilde \nu} \big( G_{\overline \chi} ^2 / N \big)
  \sum_{\substack{
        \beta (N)^\times \\
        \gamma (N)^\times
       }}
  \overline \chi(\beta \gamma^*)
  \chi (\gamma) \,
  e\big( (\beta \gamma^* + {\tilde \nu}) \gamma/N \big)
\displaybreak[0]
\\[6pt]
&=
  \chi(-m) \, {\tilde \nu} \big( G_{\overline \chi} ^2 / N \big)
  \sum_{\substack{
        \beta (N)^\times \\
        \gamma (N)^\times
       }}
  \overline \chi(\beta)
  \chi (\gamma)^2 \,
  e\big( (\beta + {\tilde \nu} \gamma)/N \big)
\displaybreak[0]
\\[6pt]
&=
  \chi(-m) \, {\tilde \nu} \big( G_{\overline \chi} ^3 / N \big)
  \sum_{\gamma (N)^\times}
  \chi (\gamma)^2 \,
  e\big( {\tilde \nu} \gamma/N \big)
\displaybreak[0]
\\[6pt]
&=
  \chi(-m) \, {\tilde \nu} \big( G_{\overline \chi} ^3 / N \big)
  \sum_{\substack{
        \gamma (L R)^\times \\
        r (N/L R)
       }}
  \big( \mathbbm{1}_R \chi^{(2)} \big) (\gamma) \,
  e\big( {\tilde \nu} (r L R + \gamma)/N \big)
\text{.}
\end{align*}
The sum ranging over $r$ vanishes, if ${\tilde \nu} \ne \nu$.  Otherwise, we have
\begin{align*}
 A_{\chi, \nu} (m)
&=
  \chi(-m) \, {\tilde \nu} \big( G_{\overline \chi} ^3  / L R \big)
  \sum_{\gamma (L R)^\times}
  \big( \mathbbm{1}_R \chi^{(2)} \big) (\gamma)
  e\big( \nu \gamma/N \big)
\\[6pt]
&=
  \chi(-m) \, \chi^{(2)}(R) \, \overline {\chi^{(2)}} \big(L R \nu / N \big) \,
  {\tilde \nu} \big( G_{\overline \chi} ^3 / L R \big) G_{\chi^{(2)}} \,
  \mu(r) \varphi\big( R / r \big)
\text{,}
\end{align*}
where
\[
  r
=
  \frac{R}{\gcd(R, L R \nu /N)}
=
  \frac{N}{\gcd(R, L R \nu / N) \, N / R}
=
  \frac{N}{\gcd(N, L \nu)}
\text{.}
\]
\end{proof}

Combining Lemma \ref{sndintreprfourierjacobiseries} and \ref{Achisum}, we deduce the next theorem.
\begin{theorem}
\label{sndintregralrepresentation}
We have
\begin{align*}
&
  \Big\{
   \mathbb{E} ^{(p) *}_{N p / \nu ,\, N p / \nu,\, \overline \chi^2}
   (\cdot,s) \,
   f_\chi \big| W_{N p}, \,
   g\big|W_{N p}
  \Big\}
\\[6pt]
&\quad=
  \frac{2}{i_N N^2}
  \big( \pi / N^2 \big)^{k - 2} \,
  p^{3s/2 - 1} \big( 1 + p^{-s} \big) \,
  {\overline \chi}^{(2)}(R) \, \mu(R) \,
  G_{\chi} ^4 G_{{\overline \chi}^{(2)}} \,
  \mathbb{D}_{f,g,\overline\chi} (s + k - 2)
\text{.}
\end{align*}
\end{theorem}
\begin{proof}
The integral is well-defined and converges absolutely due to the polynomial growth of the Eisenstein series.  Thus, using (ii) of Lemma \ref{la:eisensteinreductionproperties} and the definition of the Petersson scalar product, we find
\begin{align}
\label{eq:sndintegralrepresentation_scalarproduct}
  H(s)
&:=
  \Big\{
   \mathbb{E} ^{(p) *}_{N p / \nu,\, N p / \nu,\, \overline \chi^2}
   (\cdot,s) \,
   f_\chi \big| W_{N p}, \,
   g\big|W_{N p}
  \Big\}
\\[6pt]\nonumber
&\hphantom{:}=
  \Big\{
   \mathbb{E} ^{(p) *}_{N p,\, N \nu p,\, \overline \chi^2}
   (\cdot,s), \,
   \mathrm{trace} \Big(
   \overline {f_\chi \big| W_{N \nu p} } \;
   g \big| W_{N \nu p} \;
   \big(\det \Im(\cdot)\big)^k
   \Big)
  \Big\}
\text{.}
\end{align}
The trace operator is
\begin{align*}
  \mathrm{trace} \,:\,
&
   [\Gamma_{2,1} ^{(p) *} (N \nu p,\,
                           (N \nu)^2 p),\,
    0,\,\big(\overline \chi ^2\big)^+]
  \longrightarrow
   [\Gamma_{2,1} ^{(p) *} (N p,\,
                           N \nu p),\,
    0,\, \big(\overline \chi ^2\big)^+],
\\&
  f
  \mapsto
  (N \nu^3)^{-1} \hspace{-2em}
  \sum_{M : \Gamma_{2,1} ^{(p) *} (N \nu p,\,
                                   (N \nu)^2 p)
            \backslash
            \Gamma_{2,1} ^{(p) *} (N p,\,
                                   N \nu p)}
  \hspace{-3em}
  \big(\overline \chi ^2\big)^+ (M) \, f \big| M
\text{.}
\end{align*}
Here,
\begin{gather*}
  N \nu^3
=
  \big[ \Gamma_{2,1} ^{(p) *} (N p,\, N \nu p)
   \,:\,
        \Gamma_{2,1} ^{(p) *} (N \nu p,\, (N \nu)^2 p)
  \big]
\end{gather*}
is the index of the latter in the former group.

Set $\theta = N / \nu$ and recall the definition of $M_{d, \gamma}$, given in Section~\ref{sec:paramodulargroup}.  Define $S$ by
\begin{gather*}
  (N \nu^3)^{-1} (\det \Im(\cdot))^k \,
  \ov{S}
=
  \mathrm{trace} \Big(
   \overline {f_\chi \big| W_{N \nu p} } \,
   g \big| W_{N \nu p} \,
   (\det \Im(\cdot))^k
  \Big)
\text{.}
\end{gather*}
We write $\{1,\ldots,\nu / d\}^{2\times}$ for the set of all coprime pairs with entries between $1$ and $\nu / d$.  By sums ranging over $B \isdiv N \nu$ we mean the sum over the positive divisors of $N \nu$.  With this notation at hand, we find
\begin{align}
\label{eq:GchiS_equation}
  G_{\overline \chi}\, S
&=
  \sum_{\substack{
        d \isdiv \nu \\
        \lambda \in \{1, \ldots , \nu / d \}^{2 \times} \\
        \beta (N \nu) \\
        \mu (N)
       }}
  \hspace{-1em}
  \overline \chi (\mu) \,
  f \big| W_{N \nu p}
           M_{d,\beta - \nu p \mu} M_\lambda \;
  \overline{
   g \big| W_{N \nu p}
            M_{d,\beta} M_\lambda }
\displaybreak[0]
\\[6pt]\nonumber
&=
  \sum_{\substack{
        d \isdiv \nu \\
        \lambda \in \{1, \ldots , \nu / d \}^{2 \times} \\
        B \isdiv N \nu,\,
        \beta (N \nu / B)^\times \\
        C \isdiv N \nu,\,
        \gamma (N \nu / C)^\times \\
        B \beta \equiv C \gamma (\nu)
       }}
  \hspace{-1.5em}
  \overline \chi \Big(\frac{B \beta - C \gamma}{\nu}\Big) \,
  f \big| W_{N \nu p}
           M_{d,C \gamma} M_\lambda \;
  \overline{
   g \big| W_{N \nu p}
            M_{d,B \beta} M_\lambda }
\text{.}
\end{align}

With these computations in mind, we reconsider \eqref{eq:sndintegralrepresentation_scalarproduct}.  Apply the unfolding trick, and notice that $S$ is invariant under $D_\epsilon$ for $\epsilon = \pm 1$.  Hence we may double the integral, summing over the images of the integrand under $D_{\pm 1}$.
\begin{align}
\label{eq:H_expanded_integral}
&
  \Big(\frac{1}{2 i_N \nu^4} (N^2 / \pi)^s \, p^{3s/2} (1 + p^{-s}) \,
        \Gamma(s) L(2s, \chi^2) \, G_{\overline \chi}^{-1}
  \Big)^{-1}
  H(s)
\\[6pt] \nonumber
&\qquad=
  \int_{\cF^{(p)}_{2,1}}
  y^{k - 3} \big(y' - y^{-1} v^2 \big)^{s + k - 3}
  \sum_{\epsilon \in \{\pm 1\}}
  G_{\overline \chi} S \big| D_\epsilon (Z)
  \; dX\, dY
\text{.}
\end{align}
We will show that certain parts of this integral vanish.  More precisely, we will show that only the addends with $B = C = 1$ contribute to \eqref{eq:H_expanded_integral}.  Since $S$ is a finite sum of cusp forms, the integral \eqref{eq:H_expanded_integral} over the Fourier series converges absolutely and we may reorder it.

Expand the integrand using Lemma \ref{sndintreprfourierjacobiseries}.  Apply this Lemma to $f$ and add the superscript~$(B)$ to $\mu$ and the subscript~$B$ to $m$.  Further, write $\beta$ and $\beta^*$ instead of $\gamma$ and $\gamma^*$.  Apply the Lemma also to $g$ and add the superscript~$(C)$ and the subscript~$C$ to $\mu$ and $m$.  Whenever $\nu \ne 4$, we choose $\epsilon_\beta, \epsilon_\gamma = 1$.  Recall that, otherwise, we have $\epsilon_\beta, \epsilon_\gamma = \pm 1$.

Combining \eqref{eq:GchiS_equation} and \eqref{eq:H_expanded_integral} and inserting the Fourier-Jacobi expansions of the transforms of $f$ and $g$, the following inner integral occurs:
\begin{align}
\label{eq:H_integraloverx}
  \int_{x \in [0,\frac{1}{p}]}
  e\Big( x
     \big(
      (m_C
      C^2 / {{\mu ^{(C)}} ^2})
      -
      (m_B
      B^2 / {{\mu ^{(B)}} ^2})
     \big) \Big)
\text{.}
\end{align}
By Lemma~\ref{sndintreprfourierjacobiseries}, we may assume that $p \isdiv m_B, m_C$.  Hence
\begin{gather}
\label{sndintegralmBnCequation}
  m_C\,
  C^2 / {{\mu ^{(C)}} ^2}
=
  m_B\,
  B^2 / {{\mu ^{(B)}} ^2}
\text{,}
\end{gather}
whenever the integral~\eqref{eq:H_integraloverx} does not vanish.

We show that the contribution to \eqref{eq:H_expanded_integral} vanishes, if $2 \mu^{(B)} \nu \isdiv B$ or $2 \mu^{(C)} \nu \isdiv C$.  For reasons of symmetry, we may assume the latter.  Since $2 \nu \isdiv C$ and $C \isdiv N \nu$, we find $2 \isdiv N$.  Consequently, we have $2 \nu \isdiv N$, because $\chi$ is primitive.  Since $\beta B \equiv \gamma C\; (\nu)$, we find that $\nu \isdiv B$.  We may assume that the exponential valuations $\nu_2(B)$ and $\nu_2 (\nu)$ at $2$ are equal, since otherwise the term $\overline\chi \big( (B \beta - C \gamma)/\nu \big)$, showing up in~\eqref{eq:GchiS_equation}, vanishes.  Since $2 \nu \isdiv N$, we have $B \isdiv N$, and we conclude $\mu^{(B)} = 1$.

From \eqref{sndintegralmBnCequation}, we deduce $4 \isdiv m_B$ in case \eqref{eq:H_integraloverx} does not vanish.  We consider the sum over $\beta$ in \eqref{eq:GchiS_equation}.  Given $\beta$, choose $\widetilde{\beta}$ in the system of representatives that the sum ranges over such that $\widetilde{\beta} \equiv \beta - N \nu/2 B\;(N \nu/B)$.  From this condition, it follows that $\gcd(\widetilde{\beta}, N \nu / B) = 1$ and $\epsilon_\beta = \epsilon_{\widetilde \beta}$.  Since $\beta$ and $\beta^*$ are odd, a direct calculation yields
\begin{gather*}
  {\widetilde \beta}^*
\equiv
  \beta^* + N \nu / 2 B \quad \big(N \nu / B \big)
\end{gather*}
for any ${\widetilde \beta}^*$ satisfying $\widetilde{\beta} {\widetilde \beta}^* \equiv 1 \big(N \nu / B \big)$.  Furthermore, since $4 \isdiv m_B$, we have
\begin{equation*}
  \frac{B {\mu^{(B) *}}^2 N \nu/2 B}
       {N \nu}
  m_B \in \ZZ
\text{.}
\end{equation*}
That is,
\begin{gather*}
  e\big( {\mu^{(B) *}}^2 \beta^* p^* m_B / N \nu \big)
=
  e\big( {\mu^{(B) *}}^2 \widetilde{\beta}^* p^* m_B / N \nu \big)
\text{,}
\end{gather*}
which is the exponential term occurring when applying Lemma~\ref{sndintreprfourierjacobiseries} to \eqref{eq:GchiS_equation}.  Since ${\ov \chi} \big( (B \beta - C \gamma) / \nu \big) = - {\ov \chi} \big( (B \widetilde{\beta} - C \gamma) / \nu \big)$, the sum over all $\beta$ vanishes.

Summarizing this part of the proof, we may always assume that $2 \mu^{(C)} \nu \nisdiv C$ and $2 \mu^{(B)} \nu \nisdiv B$.  In particular, we have $\mu ^{(B)} = \mu ^{(C)} = 1$:  Indeed, assume that $\mu ^{(B)} \ne 1$ or $\mu ^{(C)} \ne 1$.  For reasons of symmetry, we may suppose that $\mu ^{(C)} \ne 1$.  Since $\nu \ne 1$, it follows from $\mu ^{(C)} \ne 1$ that $2 \nu \isdiv N$.  From the equality $\nu_2 (C) = \nu_2 ( N d p \mu^{(C)})$, it follows that $2 \mu^{(C)} \nu \isdiv C$.

Set
\begin{gather*}
  \delta = (B, C),\; b = B / \delta,\; c = C / \delta
\text{.}
\end{gather*}
If \eqref{eq:H_integraloverx} does not vanish, we have $m_B = c^2 m$ and $m_C = b^2 m$ for some $m \in \NN$.  We may assume that $\delta \isdiv \nu$, since otherwise $\overline\chi \big( B (\beta - C \gamma) / \nu \big)$ vanishes.  We split off the following inner sum of \eqref{eq:GchiS_equation}:
\begin{multline}
\label{eq:Hinnersumdecomposition}
   \sum_{\substack{
         \beta' (N \nu / \delta b c)^\times \\
         \gamma' (N \nu / \delta b c)^\times
        }}
   \sum_{\substack{
         \beta (N \nu / \delta b)^\times \\
         \gamma (N \nu / \delta c)^\times \\
                \beta - \beta'
         \equiv \gamma - \gamma'
         \equiv 0\; (N \nu / \delta b c) \\
         \delta b \beta \equiv \delta c \gamma\; (\nu)
        }}
   \hspace{-1.5em}
   \ov{\chi}\big((\delta b \beta - \delta c \gamma) / \nu\big)
\\[4pt]
   \cdot
   f_{b^2 m, \delta c, d, \lambda, \epsilon, \epsilon_{\gamma}} \,
   g_{c^2 m, \delta b, d, \lambda, \epsilon, \epsilon_{\beta}} \,
   e\big(m \delta b c (b \gamma^* - c \beta^*) p^* / (N \nu) \big)
\text{.}
\end{multline}
Observe that 
\begin{gather*}
  \delta b c
=
  \prod_{q \isdiv B C\text{,} \;q\text{ prime}}
  \hspace{-0.75em}
  q^{\max\{\nu_p(B),\,\nu_p(C)\}}
\text{.}
\end{gather*}
We conclude that $\nu \isdiv N \nu/\delta b c$. Thus, if $\beta \equiv \beta' \; (N \nu/\delta b c )$, we have $\epsilon_\beta = \epsilon_{\beta'}$.  Analogously, if $\gamma \equiv \gamma' \; (N \nu/\delta b c )$, we have $\epsilon_\gamma = \epsilon_{\gamma'}$.  Moreover, it suffices to impose the constraint $\delta b \beta \equiv \delta c \gamma\; (\nu)$ on $\beta'$ and $\gamma'$ only.

We treat the character sum in \eqref{eq:Hinnersumdecomposition}.  In analogy with $\beta^*$ and $\gamma^*$, choose ${\beta'}^*$ and ${\gamma'}^*$ such that $\beta' {\beta'}^* \equiv 1 \; (N \nu / \delta b)$ and $\gamma' \gamma'^* \equiv 1 \; (N \nu / \delta c)$.  Note that \mbox{$\beta \equiv \beta' \; (N \nu/\delta b c)$} if and only if $\beta^* \equiv \beta'^* \; (N \nu/\delta b c)$.  Analogously, $\gamma \equiv \gamma' \; (N \nu / \delta b c)$ if and only if $\gamma^* \equiv \gamma'^* \; (N \nu/\delta b c)$.  For fixed $\beta'$ and $\gamma'$, the exponents in \eqref{eq:Hinnersumdecomposition} differ by multiples of
\begin{gather*}
  \frac{m \, \delta b c \, b \, (N \nu / \delta b c)}{N \nu}
=
  bm \in \ZZ,
\qquad
  \frac{m \, \delta b c \, c \, (N \nu / \delta b c)}{N \nu}
=
  cm \in \ZZ
\text{.}
\end{gather*}
We may assume that $\delta b \beta' \equiv \delta c \gamma' \;(\nu)$.  We focus on the following inner sum of \eqref{eq:Hinnersumdecomposition}:
\begin{align}
\label{eq:Hinnercharactersum}
\sum_{\substack{
        \beta (N \nu / \delta b)^\times \\
        \gamma (N \nu / \delta c)^\times \\
               \beta - \beta'
        \equiv \gamma - \gamma'
        \equiv 0\; (N \nu / \delta b c)
       }}
\hspace{-1.5em}
  \overline \chi \big( (\delta b \beta - \delta c \gamma / \nu \big)
\text{.}
\end{align}
If $\delta b c \nisdiv \nu$, then $N \nisdiv N \nu / \delta b c$ and \eqref{eq:Hinnercharactersum} vanishes.  Thus we may assume that $\delta b c \isdiv \nu$.  Since $\delta b \beta \equiv \delta c \gamma\; (\nu)$, we have $\gcd(c, \nu) = \gcd(b, \nu) = 1$.  Hence $b = c = 1$ and $\delta \isdiv \nu$.  If $\nu \ne 1$, the expression $\chi\big( (\beta - \gamma) / (\nu / \delta) \big)$ vanishes, whenever $\delta = \nu$, because $\beta$ and $\gamma$ are odd.  Thus if $\nu \ne 1$, we may assume that $\delta \ne \nu$.

With the conclusions drawn so far in mind, we rewrite another part of \eqref{eq:GchiS_equation}.  In \eqref{eq:Hinnersumdecomposition}, observe that $\epsilon_\beta = \epsilon_{\beta'}$, whenever $\beta \equiv \beta' \;(\nu)$.  Analogously, $\epsilon_\gamma = \epsilon_{\gamma'}$, whenever $\gamma \equiv \gamma' \;(\nu)$.  In addition, we find $\epsilon_\gamma = \epsilon_\beta$ in all cases but $\nu = 4$ and $\delta = 2$.  In this case, we may assume that $\epsilon_\gamma \ne \epsilon_\beta$, since otherwise $\chi \big( (b \beta - c \gamma) / (\nu / \delta) \big) = 0$.  We consider the following inner sum of \eqref{eq:GchiS_equation}:
\begin{align}
\label{eq:finalcharactersum}
&
  \sum_{\substack{
         \epsilon \in \{\pm 1\} \\
         \beta (N \nu / \delta)^\times \\
         \gamma (N \nu / \delta)^\times \\
         \beta \equiv \gamma \; (\nu / \delta) \\
        }}
    \hspace{-0.5em}
  \overline \chi \big( (\beta - \gamma) (\nu / \delta) \big) \,
  f_{m, \delta , d, \lambda, \epsilon, \epsilon_{\gamma'}} (\tau, z) \,
  \overline{g_{m, \delta, d, \lambda, \epsilon, \epsilon_{\beta'}} (\tau, z)} \,
\\[4pt]\nonumber
&\hspace{5em}
  \cdot e\big( m (\gamma^* - \beta^*) p^* / (N \nu / \delta) \big)
\text{.}
\end{align}
Recall that we assume that $C \nu, B \nu = \delta \nu \isdiv N$.  By Lemma \ref{sndintreprfourierjacobiseries}, $f_{m, \delta , d, \lambda, \epsilon, \epsilon_{\gamma'}}$ and $g_{m, \delta, d, \lambda, \epsilon, \epsilon_{\beta'}}$  only depend on $\epsilon,\,\epsilon_\gamma$ and $\epsilon_\beta$ by means of the products $\epsilon \epsilon_\gamma$ and $\epsilon \epsilon_\beta$.  Hence we may reorder the sum and split it such that it can be evaluated using Lemma \ref{Achisum}.  By this Lemma, \eqref{eq:finalcharactersum} vanishes, if $\delta \ne 1$.  Otherwise, it equals
\begin{align*}
& \quad
  \sum_{\epsilon \in \{\pm 1\}}
  A_{\overline \chi, \nu}(m p^*) \,
  f_{m, \delta , d, \lambda, 1, \epsilon} (\tau, z) \,
  \overline{g_{m, \delta, d, \lambda, 1, \epsilon} (\tau, z)}
\text{.}
\end{align*}

So far, we have proved that any contribution to \eqref{eq:H_expanded_integral} vanishes, unless \mbox{$B = C = 1$}.  By Lemma \ref{sndintreprfourierjacobiseries}, explicit expressions for the remaining $f_{m, \delta, d, \lambda, 1, \epsilon}$ and $g_{m, \delta, d, \lambda, 1, \epsilon}$ are known.  They depend on $d$ and $\lambda$ by means of $d \lambda_1$ and $d \lambda_2$ only.

Before simplifying the remaining terms in \eqref{eq:H_expanded_integral}, set $\Lambda = \{1,\ldots, \nu\}^{2}$.  We will make use of the periodicity of $f_m$ and $g_m$ in the elliptic variable $w$.
\begin{align*}
& \quad
  \Big(
   \frac{1}{2 i_N \nu^4} \big(N^2 / \pi \big)^s \, p^{3s/2} (1 + p^{-s}) \,
   \Gamma(s) L(2s, \chi^2) \,
   G_{\overline \chi}^{-1}
  \Big)^{-1}
  H(s)
\\[6pt]
&=
  \int_{\cF^{(p)}_{2,1}}
  y^{k - 3} \big(y' - y^{-1} v^2\big)^{s + k - 3}
  \sum_{\epsilon \in \{\pm 1\}}
  G_{\overline \chi} S \big| D_\epsilon (Z) \;dX\, dY
\allowdisplaybreaks
\\[6pt]
&=
  (4 \pi)^{-s - k + 2} \, p^{-1} \, \Gamma(s + k - 2)
  \int_{\cF^*}
  y^{k - 3}
  \sum_{\substack{
        m \ge 1 \\[1pt]
        \epsilon \in \{\pm 1\},\, \lambda \in \Lambda
       }}
  A_{\overline \chi, \nu}(m) \,
  m^{- s - k + 2} \,\cdot
\\[4pt]
& \qquad\quad
  \cdot
  f_m \Big(
       \tau,
       z
       + (p p^* \epsilon / \nu)  \big( \tau \lambda_1 + \lambda_2 \big)
      \Big)
  \,
  \overline{
   g_m \Big(
        \tau,
        z
        + (p p^* \epsilon / \nu) \big( \tau \lambda_1 + \lambda_2 \big)
       \Big)
  }
\\[4pt]
&\qquad\quad
  \cdot
  e^{- 4 \pi m \big(  (p p^* \epsilon / \nu )^2 y \lambda_1 ^2
                    + 2 p p^* \epsilon  v \lambda_1 / \nu + y^{-1} v^2 \big)}
  \;dx\, dy\, du\, dv
\\[6pt]
&=
 (4 \pi)^{-s - k + 2} \,p^{-1} \, \Gamma(s + k - 2)
  \int_{\cF^*}
  y^{k - 3}
  \sum_{\substack{
        m \ge 1 \\[1pt]
        \epsilon \in \{\pm 1\},\, \lambda \in \Lambda
       }}
  A_{\overline \chi, \nu}(m) \,
  m^{- s - k + 2} \,\cdot
\\[4pt]
& \qquad\quad
  \cdot
  f_m (\tau, z) \,
  \overline{
   g_m (\tau, z)
  } \,
  e^{- 4 \pi m y^{-1} v^2}
  \;dx\, dy\, du\, dv
\\
&=
  4 \nu^2 A_{\overline \chi, \nu}(1) \,
  (4 \pi)^{-s - k + 2} \, p^{-1} \, \Gamma(s + k - 2) \,
  D_{f,g, \overline \chi} (s + k - 2)
\text{.}
\end{align*}

Using  Lemma \ref{Achisum}, we find that $H(s)$ equals
\begin{align*}
&\quad
  \frac{1}{2 i_N \nu^4}
  \big( N^2 / \pi \big)^s \, p^{3s/2} (1 + p^{-s}) \,
  \Gamma(s) L(2s, \chi^2) \, G_{\overline \chi}^{-1} \,\cdot
\\[4pt]
&\qquad\quad
  \cdot
  4 \nu^2 A_{\overline \chi, \nu} (1) \,
  (4 \pi)^{-s - k + 2} \,  p^{-1} \, \Gamma(s + k - 2) \,
  D_{f,g, \overline \chi} (s + k - 2)
\\[6pt]
&=
  \frac{2}{i_N N^2}
  \big( \pi / N^2 \big)^{k - 2} \,
  p^{3s/2 - 1} (1 + p^{-s}) \,
  {\overline \chi}^{(2)}(R) \, \mu(R) \,
  G_{\chi}^4 G_{{\overline \chi}^{(2)}} \,
  \mathbb{D}_{f,g, \overline \chi} (s + k - 2)
\text{.}
\end{align*}
\end{proof}


\section{The Rankin convolution and the spinor $L$-function}
\label{sec:spinorLseries}

The functional equation of the twisted Rankin convolution is a direct consequence of the results in Section \ref{sec:fstintegralrepresentation} and \ref{sec:sndintegralrepresentation}.  From the next corollary, we will deduce the functional equation of the twisted spinor $L$\nbd function.
\begin{corollary}
\label{functionaleqnconvolution}
Suppose that $f, g \in \big[\Gamma ^{(p) *}, k, \mathbbm{1}_1 ^{k-}\big]_0$ and adopt the assumptions on $\chi$ and $N$ from Section \ref{sec:sndintegralrepresentation}.  Then we have
\begin{equation*}
  \mathbb{D}_{f,g,\chi} (s)
=
  \frac{G_\chi ^4}{N^2} \,
  p^{3 (k - s - 1)}
  (1 + p^{-k + s})
  (1 + p^{-s + k - 2})^{-1}
  \mathbb{D}_{f,g,\overline \chi} (2k - 2 - s)
\text{.}
\end{equation*}
\end{corollary}

The rest of this section will be dedicated to studying the twisted spinor $L$-series and its completion.  Given a Siegel eigenform $f$ denote the eigenvalues under the Hecke operators $T(p)$ and $T(p^2)$ by $\lambda_f(p)$ and $\lambda_f(p^2)$ (see \citep{Kr90} for a definition of the Hecke operators).
\begin{definition}
\label{def:spinorLdefinition}
The spinor $L$-series attached to a weight $k$ cuspidal Siegel eigenform for $\Sp{2}(\ZZ)$ twisted by a Dirichlet character $N$ is
\begin{multline}
  Z_f^\chi (s)
:=
  \prod_p
  \Big(
     1
   - \lambda_f(p) \, \chi(p) p^{-s}
   + \big(\lambda_f(p)^2 - \lambda_f(p^2) - p^{2 k - 4}\big) \big(\chi(p) p^{-s}\big)^2
\\[4pt]
   - \lambda_f(p) p^{2 k - 3} \big(\chi(p) p^{-s}\big)^3
   + p^{4 k - 6} \big(\chi(p) p^{-s}\big)^4
  \Big)
\text{.}
\end{multline}
Denote the period of $\chi$ by $N$.  We call
\begin{gather}
  \bbZ_f^\chi(s)
:=
  (2 \pi / N)^{-2 s} \, \Gamma(s) \Gamma(s - k + 2) \, Z_f^\chi
\end{gather}
the completed twisted spinor $L$-series attached to $f$.
\end{definition}

We relate the spinor $L$-function to the above Rankin convolution, using a representation given by Gritsenko in \citep{gritsenkoarithmeticallifting}.  Unfortunately, he did not make this representation explicit.  We indicate how to obtain it, adopting his notation.  Assume that $f$ is a holomorphic Siegel eigenform that is a cusp form such that the Fourier-Jacobi coefficients satisfy $f_1 = 0$ and $f_p \ne 0$.  By \citespinorgritsenko, we have
\begin{align*}
& \quad
  L(2 s - 2 k + 4, \mathbbm{1}_p) \, D_{\mathrm{Sym}^p f, F_{f_p}}
\displaybreak[0]
\\[6pt]
&=
  L(2 s - 2 k + 4, \mathbbm{1}_p) \,
  \sum_{m \ge 1}
  \big\langle f_{m p}^*,\, m^{2-k} f_p \big|_k T_{-} ^{(1)} (m) \big\rangle
  (m p)^{-s}
\displaybreak[0]
\\[6pt]
&=
  p^{-s} L(2 s - 2 k + 4, \mathbbm{1}_p) \,
  \sum_{m \ge 1}
  \big\langle f_{m p} ^* \big|_k T_{+} ^{(1)} (m),\, f_p \big\rangle
  m^{-s}
\displaybreak[0]
\\[6pt]
&=
  p^{-s}
  \Big(
  (p - p^{2 k - 2 s - 1}) \langle f_p, f_p \rangle
\\[4pt]
& \hspace{5em}
  +
  \big((-1)^k p^{3 - s - k} + p^{2 - 2 s}\big)
  \big\langle f_p \big|_k T_{-} ^{(1)} (p) T_{+} ^{(1)} (p),\, f_p \big\rangle
  \Big) \,
  Z_f (s)
\\[6pt]
&=
  p^{-s}
  (1 - p^{k - 2 - s})(p + (-1)^k p^{k - s}) \,
  \langle f_p, f_p \rangle \,
  Z_f (s)
\text{.}
\end{align*}

For the time being, we assume that $p \equiv 1 \,(N)$ and $f_p \ne 0$.  Starting with Corollary~\ref{functionaleqnconvolution}, a straightforward calculation yields
\begin{align*}
  \ZZ_f ^{\chi} (s)
=
  (-1)^k
  \frac{G_{\chi} ^4}{N^2} \,
  \ZZ_f ^{\overline \chi} (2 k - 2 - s)
\text{.}
\end{align*}

To prove the Main Theorem we are left with considering the Fourier\nbd Jacobi expansion of an arbitrary Siegel eigenform.
\begin{proof}[{\tit Proof of the Main Theorem}]
The functional equation was proved in \citep{kussthesis}, if \mbox{$f_1 \ne 0$}.  Hence we may suppose that $f_1 = 0$.  By the preceding calculations, it suffices to prove that for each $N \in \NN$ there is a prime $p_N \equiv 1\;(N)$ such that $f_{p_N} \ne 0$.

Fix such an $N$.  We will write
\begin{gather}
  f(Z)
=
  \sum_{0 < T \in (\frac{1}{2}\ZZ)^{(2,2)}} a(T) \, e^{2 \pi i\, \tr(T Z)}
\end{gather}
for the Fourier expansion of $f$.

We first show that there is a non\nbd vanishing Fourier coefficient $a(T)$ of $f$ such that $T$ represents a unit and square modulo $N$.  Second, we will show that we may choose a primitive such $T$.  Finally, we will use the correspondence between values modulo $N$ represented by a quadratic form and ray classes of a quadratic number field to prove the claim.

Suppose that there is no Fourier index $T$ with non\nbd vanishing coefficient $a(T)$ attached to it that represents a unit and square modulo $N$.  We consider the Taylor expansion at $z = 0$ of an arbitrary Fourier-Jacobi coefficient $f_m$.  By our assumption $f_m(\tau, 0)$ is lacunary and has weight greater than $4$.  That is, $f_m(\tau, 0) = 0$.  Assume that we have proved $\partial_z^l f_m (\tau, 0) = 0$ for all $0 \le l < L$ with $L \in \NN$ given.  Then $\partial_z^L f_m(\tau, 0)$ is a lacunary modular form for $\SL{2}(\ZZ)$ of weight greater than $4$.  We conclude that $\partial_z^L f_m(\tau, 0) = 0$.  By mathematical induction, it follows that all derivatives of $f_m$ with respect to $z$ vanish at $z = 0$, and thus $f_m = 0$.  Since $m$ was arbitrary, $f = 0$, contradicting the assumptions.

We have proved that there is at least one index $T$ with $a(T) \ne 0$ that represents a unit and square modulo $N$.  We claim that there is such an index $T$ that is primitive.  Let $r$ be the minimal square and a unit modulo $N$ represented by some, non\nbd fixed $T$ satisfying $a(T) \ne 0$.  Consider the proof of \cite[Theorem 4.7]{gritsenkoparamodularforms}.  Adopting Gritsenko's notation, we find that the identities
\begin{align*}
  f_d \{F\big|_k T(e) \} (\tau, z)
&=
 \big( f_r \{F\}\big|_k T_+ (e) \big) (\tau, z)
\quad\text{and}
\\
  f_r(\tau, z)\, \exp(2 \pi i\, r \omega)
&=
  \sum_{\substack{N = \left(\begin{smallmatrix}* & * \\ * & r\end{smallmatrix}\right) \in \mathfrak{B}_2}} \!\!
  a(N) \exp(2 \pi i\, \tr(N Z))
\end{align*}
only involve quadratic forms that represent a unit and square modulo $N$.  Indeed,  $(r / e) (a / b)$ is a unit and square modulo $N$, if $a b = e \isdiv r$.  Using this fact, the calculations carried out by Gritsenko prove the existence of a primitive $T$ representing a unit and square modulo $N$, as we have claimed.

The index $T$ also represents an $a_N \equiv 1\;(N)$.  Indeed, it corresponds to an ideal class $\mathfrak{k}$ in the narrow sense of the imaginary quadratic field $K = \QQ\big(\sqrt{-\det(T)}\big)$.  To prove the existence of $a_N$ it suffices to prove that there is an element of $\mathfrak{k}$ that has norm $a_N$ over $\QQ$.  We consider the ray class group associated to the principal ideal $(N) \subseteq \mathfrak{O}_K$ in the ring of integers of $K$. We know that there is a ray class $\mathfrak{k}'$ containing an element that has norm $a_N \equiv 1\;(N)$.  By \citespinorneukirch and an argument paralleling the one in \citespinorzagier, there are infinitely many split prime ideals in $\mathfrak{k}'$.  This proves the existence of $p_N$.
\end{proof}


\bibliographystyle{amsalpha}

\bibliography{bibliography}

\end{document}